\documentclass[11pt]{article}
\usepackage[margin=1in]{geometry}
\usepackage[T1]{fontenc}
\usepackage{lmodern}
\usepackage{amsmath,amssymb,amsthm,mathtools}
\usepackage{booktabs,tabularx,array,enumitem}
\usepackage{microtype}
\usepackage{graphicx}
\usepackage{subcaption}
\usepackage{fancyhdr}
\usepackage[hidelinks]{hyperref}
\usepackage[nameinlink,noabbrev]{cleveref}
\numberwithin{equation}{section}
\newtheorem{theorem}{Theorem}[section]
\newtheorem{lemma}[theorem]{Lemma}
\newtheorem{proposition}[theorem]{Proposition}

\theoremstyle{definition}
\newtheorem{assumption}[theorem]{Assumption}
\newtheorem{remark}[theorem]{Remark}
\newcommand{\bu}{\boldsymbol u}
\newcommand{\bv}{\boldsymbol v}
\newcommand{\bw}{\boldsymbol w}
\newcommand{\bxi}{\boldsymbol\xi}
\newcommand{\be}{\boldsymbol\eta}
\newcommand{\bz}{\boldsymbol\zeta}
\newcommand{\bn}{\boldsymbol n}
\newcommand{\bD}{\boldsymbol D}
\newcommand{\bP}{\boldsymbol P}
\newcommand{\bs}{\boldsymbol\sigma}
\newcommand{\br}{\boldsymbol r}

\newcommand{\dt}{D_t}
\newcommand{\Piu}{\Pi_u}
\newcommand{\Pip}{\Pi_p}
\newcommand{\Pieta}{\Pi_\eta}
\newcommand{\Pixi}{\Pi_\xi}
\newcommand{\Piphi}{\Pi_\phi}
\newcommand{\ehv}[2]{\boldsymbol e_{#1}^{h\ifx\relax#2\relax\else,#2\fi}}
\newcommand{\ehs}[2]{e_{#1}^{h\ifx\relax#2\relax\else,#2\fi}}
\newcommand{\eIv}[2]{\boldsymbol e_{#1}^{I\ifx\relax#2\relax\else,#2\fi}}
\newcommand{\eIs}[2]{e_{#1}^{I\ifx\relax#2\relax\else,#2\fi}}

\newcommand{\dual}[2]{\left\langle #1,#2\right\rangle_\Gamma}
\newcommand{\norm}[2]{\left\|#1\right\|_{#2}}
\newcommand{\E}{\mathcal E}
\newcommand{\N}{\mathcal N}
\newcommand{\Diss}{\mathcal D}
\newcommand{\R}{\mathcal R}
\newcommand{\Vf}{\boldsymbol V_f}
\newcommand{\Vp}{\boldsymbol V_p}
\newcommand{\Vfh}{\boldsymbol V_{f,h}}
\newcommand{\Vph}{\boldsymbol V_{p,h}}

\setlist[itemize]{leftmargin=1.4em,itemsep=3pt,topsep=4pt}
\setlist[enumerate]{leftmargin=1.6em,itemsep=3pt,topsep=4pt}
\allowdisplaybreaks[2]
\title{\fontsize{17}{21}\selectfont A Second-Order Parallel Neumann–Robin Splitting for the Stokes–Biot System with a Linear Time-Step Restriction}
\author{Yifan Wang\thanks{Department of Mathematics and Statistics, Texas Tech University, Lubbock, TX, USA. \texttt{yifan.wang@ttu.edu}.}
\and Jeonghun J.  Lee\thanks{Corresponding author. Department of Mathematics, Baylor University, Waco, TX, USA. \texttt{Jeonghun\_Lee@baylor.edu}.}
}
\date{}

\DeclareRobustCommand{\DIFaddtex}[1]{\ifmmode\expandafter\DIFaddmath\else\expandafter\DIFaddtext\fi{#1}}
\newcommand{\DIFaddmath}[1]{\begingroup\color{red}#1\endgroup}
\newcommand{\DIFaddtext}[1]{\ifvmode\leavevmode\fi{\color{red}#1}}

\providecommand{\DIFdel}[1]{}
\providecommand{\DIFaddbegin}{}
\providecommand{\DIFaddend}{}

\providecommand{\DIFdelFL}[1]{}

\begin{document}
\maketitle
\begin{abstract}
In this paper we develop a second-order fully decoupled, non-iterative partitioned method for the time-dependent Stokes--Biot problem on fixed domains. 
% The method retains the BDF2 time discretization and second-order interface extrapolation of the Robin-based scheme of Wang, Lee, and \v{C}ani\'{c} \cite{WangLeeCanic2026SecondOrder}, but replaces the poroelastic mechanical normal Robin condition by its physical Neumann-type traction condition. Pore pressure and structure velocity remain coupled implicitly within the Biot solve, while the fluid and Biot subproblems are independent and can be advanced concurrently. This modification removes the extrapolated normal structure-velocity residual without changing the continuous transmission problem or the formal temporal order. For fixed positive material and interface coefficients, we prove a fully discrete energy bound under a linear mesh-dependent restriction, $\Delta t\le c h$. The improvement over the parabolic restriction in the earlier Robin--Robin analysis follows from eliminating the interface residual that requires two inverse-trace estimates. We also establish a bulk energy-norm error bound of order $h^k+\Delta t^2$ under explicit regularity, finite element approximation, and initialization assumptions. The proof uses a common elasticity projection for displacement and velocity and a discrete summation-by-parts argument for the pressure-projection traction term. 
In previous work \cite{WangLeeCanic2026SecondOrder}, a Robin--Robin scheme with these features was introduced, and its stability and error estimates were established under the sufficient parabolic CFL-type condition $\Delta t\lesssim C_*h^2$.
% REV10-DEL: "In a previous work [1] a numerical scheme with same features was devised with Robin-Robin splitting
% approach but the numerical scheme needs a parabolic CFL-type condition Delta t <~ h^2 for the stability and error analyses."
The main improvement of the new scheme is weakening the CFL-type condition to a linear CFL-type condition $\Delta t \lesssim C_*h$. This can be achieved with a new splitting approach by replacing the poroelastic mechanical normal Robin condition by its physical Neumann-type traction condition, which results in a Neumann-Robin splitting.
For fixed coefficients satisfying the stated assumptions, we prove a fully discrete energy bound under a  restriction of $\Delta t\lesssim  C_*h$. 
The improvement over the parabolic restriction in the earlier Robin--Robin analysis follows from eliminating the interface residual that requires two inverse-trace estimates. 
We also establish  an energy-norm error bound of order $h^k+\Delta t^2$ under explicit regularity, finite element approximation, and initialization assumptions. 
Numerical experiments with manufactured solutions examine spatial and temporal accuracy and compare the Neumann–Robin and Robin–Robin schemes. A 3D cylindrical pressure-pulse benchmark further investigates their numerical stability behavior. The results distinguish the less restrictive sufficient stability condition established analytically from the observed stability thresholds of the two schemes.
\end{abstract}

\section{Introduction}
\label{sec:intro}
Fluid--poroelastic structure interaction couples free fluid motion with deformation and filtration in a saturated porous material. Such interactions arise in tissue perfusion, arterial-wall mechanics, filtration devices, and subsurface flow. A standard mathematical description combines the incompressible Stokes or Navier--Stokes equations with Biot poroelasticity \cite{Biot1941,Biot1955}. The interface conditions enforce conservation of normal fluid flux, accounting for both the porous skeleton motion and Darcy filtration, together with traction of balance and the Beavers–Joseph–Saffman tangential slip law. These conditions connect the fluid velocity and stress to both the skeleton motion and the pore-fluid response. Early computational work on the Navier--Stokes--Biot coupling includes the monolithic and partitioned formulations of Badia, Quaini, and Quarteroni \cite{BadiaQuainiQuarteroni2009}. Since then, numerical methods for Stokes-Biot problems with various formulations and finite element discretizations have been developed. We refer to \cite{Bukac2015Nitsche,BukacYotovZunino2015,AmbartsumyanEtAl2018,LiYotov2022,CaucaoLiYotov2022,RuizBaierEtAl2022,CesmeliogluLeeRhebergen2023} 
for related formulations and discretizations.

The numerical treatment depends in part on the variational formulation chosen for each component. In a \emph{primal Darcy formulation}, the filtration velocity is eliminated through Darcy's law and the pore pressure is sought in an $H^1$ space. Coupled with displacement-based elasticity, this leads to a poroelastic displacement--pressure formulation, with structure velocity introduced when inertia is retained. 
In a \emph{dual-mixed Darcy formulation}, filtration velocity and pore pressure are independent unknowns, typically in $H(\operatorname{div})$ and $L^2$, respectively. 
The normal filtration flux is then a natural trace variable, and the normal-flux transmission condition is of essential type for the mixed Darcy formulation; it may be imposed weakly, for example through a Lagrange multiplier or Nitsche-type interface operators.

Several methods have been developed for the mixed setting. Buka\v{c}, Yotov, Zakerzadeh, and Zunino \cite{Bukac2015Nitsche} used a mixed Darcy formulation and Nitsche interface operators to construct and analyze partitioned Stokes--Biot discretizations. Ambartsumyan, Khattatov, Yotov, and Zunino \cite{AmbartsumyanEtAl2018} imposed normal-flux continuity through a Lagrange multiplier and established stability and error estimates for mixed finite element approximations. Mixed formulations can also be applied to the mechanical variables: Li and Yotov \cite{LiYotov2022} combined weakly symmetric stress--displacement--rotation elasticity with mixed Darcy flow and a velocity--pressure Stokes formulation. Caucao, Li, and Yotov \cite{CaucaoLiYotov2022} developed a fully mixed Stokes--Biot formulation and a multipoint stress-flux method, using interface multipliers and local elimination of stress, rotation, and filtration variables. Thus, the terms primal and mixed must be interpreted componentwise, rather than as interchangeable descriptions of the entire coupled system.

Other formulations address complementary discretization issues. Ruiz-Baier, Taffetani, Westermeyer, and Yotov \cite{RuizBaierEtAl2022} introduced total poroelastic pressure as an additional unknown and developed a mixed-primal Biot--Stokes formulation using Stokes-stable pairs for the fluid and mechanical variables. Their transmission conditions are incorporated without interface Lagrange multipliers. Cesmelioglu, Lee, and Rhebergen \cite{CesmeliogluLeeRhebergen2023} proposed hybridizable discontinuous Galerkin methods with $H(\operatorname{div})$-conforming velocity and displacement approximations, elementwise satisfaction of the compressibility equations, and a locking-free error analysis. These developments concern the spatial formulation and interface enforcement; the design of a time-splitting method requires an additional analysis of the lagged interface exchange.

The use of Robin transmission conditions in partitioned Stokes--Biot methods is closely related to developments in incompressible fluid--structure interaction. Fern\'{a}ndez, Mullaert, and Vidrascu \cite{FernandezMullaertVidrascu2013} developed explicit Robin--Neumann schemes for incompressible fluids coupled to thin-walled elastic structures, incorporating fluid--solid inertia into the interface treatment. Burman, Durst, Fern\'{a}ndez, and Guzm\'{a}n \cite{BurmanEtAl2022} established fully discrete stability and error estimates for a non-iterative Robin--Robin fluid--structure scheme, with interfacial stresses evaluated variationally. The same authors subsequently provided a unified Robin--Robin splitting analysis for parabolic--parabolic and parabolic--hyperbolic interface problems \cite{BurmanEtAl2023}. These works provide foundations for explicit interface coupling, while their models differ from Stokes--Biot interaction, where filtration and pore pressure introduce additional transmission terms.

For fluid--poroelastic interaction itself, Buka\v{c}, Yotov, and Zunino \cite{BukacYotovZunino2015} developed a Lie operator-splitting method for a fluid interacting with a multilayered structure containing a thin elastic membrane and a thick poroelastic layer. The Nitsche-based method in \cite{Bukac2015Nitsche} instead separates the Stokes, Darcy, and elasticity equations through time lagging and stabilization of the interface operators. Dalal, Durst, Quaini, and Yotov \cite{Dalal2025Robin} developed a Robin--Robin splitting with an auxiliary interface variable representing the Robin data and one fluid and one Biot solve per time step, together with stability and time-discretization error analysis. Guo, Sun, Wang, Yue, and Zheng \cite{GuoEtAl2025} proposed a fully parallelizable loosely coupled method for moving-domain Navier--Stokes--Biot interaction, using previous-time-level interface data to advance the fluid and poroelastic subproblems independently. The partition of the Biot equations, the interface unknowns, and the order in which data are exchanged differ among these approaches.

Second-order temporal accuracy introduces further requirements on the interface approximation. Oyekole and Buka\v{c} \cite{OyekoleBukac2020} developed loosely coupled methods based on BDF2 with extrapolated coupling terms and on a Crank--Nicolson/leap-frog construction. Kunwar, Lee, and Seelman \cite{KunwarLeeSeelman2020} studied second-order time discretizations for a coupled quasi-Newtonian fluid--poroelastic system. Parrow and Buka\v{c} \cite{ParrowBukac2025} developed a strongly coupled Robin–Robin partitioned method based on a refactorized Cauchy $\theta$-like scheme. Subsequently, they analyzed a loosely coupled Robin–Robin method with conditional stability and a suboptimal temporal error estimate of $\sqrt{\Delta t}$, and investigated an iterative higher-order alternative exhibiting second-order convergence in numerical experiments \cite{ParrowBukac2026}. 

The immediate starting point of the present work is the second-order explicit splitting of Wang, Lee, and \v{C}ani\'{c} \cite{WangLeeCanic2026SecondOrder}. That method uses a primal Darcy formulation, displacement-based elasticity, BDF2 time stepping, and second-order extrapolation of the interface data. The Stokes and Biot blocks are independent at each time step and can be solved concurrently. Its mechanical normal condition on the Biot side contains an auxiliary Robin term involving the normal structure velocity. This term cancels from the continuous transmission condition, but produces a second time difference of normal structure velocity when one occurrence is imposed at the new time and the other is extrapolated. The earlier fully discrete analysis controls this contribution under a parabolic restriction, $\Delta t\lesssim C_* h^2$.

In this paper, we retain the same two-block partition and second-order time discretization, but replace the Biot mechanical normal Robin condition by the physical Neumann-type condition
\[
\boldsymbol n_p\cdot\boldsymbol\sigma_p(\boldsymbol\eta,\phi)\boldsymbol n_p+\phi=0.
\]
Both stress and pore pressure are evaluated within the current Biot solve. Since pore pressure is already a Biot unknown, this treatment preserves the independence of the fluid and poroelastic solves. The fluid normal condition, Darcy interface condition, and tangential Robin/BJS conditions retain their Robin form. We refer to the resulting algorithm as a \emph{Neumann--Robin splitting}, listing the poroelastic mechanical normal condition first. The modification concerns the normal mechanical coupling, while filtration and tangential slip remain coupled through their respective transmission conditions.

The principal analytical benefit is the removal of the extrapolated normal structure-velocity residual. In the earlier energy argument, this residual is paired with the current normal structure velocity, and estimating both traces through bulk $L^2$ norms generates the factor associated with the parabolic restriction. In the modified scheme, each remaining explicit interface term can be estimated using at most one inverse trace, together with fluid viscous dissipation, Darcy dissipation, or tangential slip control. 
For fixed positive material and interface coefficients, we establish a closed homogeneous fully discrete energy estimate under the linear restriction $\Delta t\le C_* h$. The extension to nonhomogeneous source terms is discussed separately.
Our comparison here only concerns the sufficient conditions obtained for this scheme and the specific second-order method in \cite{WangLeeCanic2026SecondOrder}, and thus it does not imply that all Robin--Robin formulations require a parabolic restriction.

We establish unique solvability of the two independent subproblems and derive an exact BDF2 discrete energy identity for the proposed Neumann–Robin algorithm. For the homogeneous problem, we prove a fully discrete energy stability estimate under a linear mesh-dependent time-step restriction. Under stated regularity, finite element approximation, and initialization assumptions, we further establish a bulk energy-norm error estimate of order \(h^k+\Delta t^2\) for the fluid velocity, structure velocity, pore pressure, and displacement. The error analysis uses a common elasticity projection for displacement and structure velocity, together with a discrete summation-by-parts argument to control the pore-pressure projection traction term without requiring normal structure-velocity dissipation. Finally, we construct a smooth polynomial manufactured solution satisfying the physical interface conditions without additional interface loads and present spatial and temporal convergence studies. These experiments assess the accuracy of the proposed splitting scheme, while the improvement in the sufficient time-step restriction is established through the theoretical stability analysis. 

% Section~\ref{sec:model} introduces the continuous problem. Section~\ref{sec:scheme} presents the Neumann--Robin algorithm and its relation to the earlier Robin--Robin formulation. Stability is established in Section~\ref{sec:stability}, and the error estimate is proved in Section~\ref{sec:error}. Section~\ref{sec:comparison} discusses the scope of the improvement, and Section~\ref{sec:verification} reports the numerical validation based on an example with manufactured solutions.
The paper is organized as follows. Section~\ref{sec:model} introduces the continuous problem  and Section~\ref{sec:scheme} presents the Neumann--Robin algorithm. Stability is established in Section~\ref{sec:stability}, and the error estimate is proved in Section~\ref{sec:error}. Section~\ref{sec:comparison} discusses the  improvement in comparison with the Robin--Robin splitting, and 
Section~\ref{sec:verification} presents spatial and temporal convergence studies using manufactured solutions, compares the Neumann-Robin and Robin-Robin splitting schemes, and investigates their numerical stability behavior in a 3D cylindrical pressure-pulse benchmark.

\section{The continuous Stokes--Biot problem}
\label{sec:model}
\subsection{Domains and equations}
Let $\Omega_f,\Omega_p\subset\mathbb R^d$, $d\in\{2,3\}$, be bounded Lipschitz domains with common interface $\Gamma$. The remaining boundaries are denoted by $\Sigma_f$ and $\Sigma_p$. Both subdomains and the interface are fixed. Their outward normals satisfy $\bn_p=-\bn_f$ on $\Gamma$, and
\[
\bP\bv:=\bv-(\bv\cdot\bn_f)\bn_f
=\bv-(\bv\cdot\bn_p)\bn_p
\]
is the common tangential projector. We assume that the exterior Dirichlet portions have positive measure and provide the Korn and Poincar\'{e} inequalities used below.

The fluid velocity and pressure are $(\bu,p_f)$, and the poroelastic displacement, velocity, and pore pressure are $(\be,\bxi,\phi)$. With $\bD(\bv)=(\nabla\bv+\nabla\bv^T)/2$, define
\begin{align}
\bs_f(\bu,p_f)&=2\mu_f\bD(\bu)-p_f\boldsymbol I,\notag\\
\bs_p(\be,\phi)&=2\mu_p\bD(\be)+\lambda_p(\nabla\cdot\be)\boldsymbol I-\alpha\phi\boldsymbol I,
\label{eq:stress}
\end{align}
where $\boldsymbol{I}$ is the $d\times d$ identity matrix. The governing equations are
\begin{equation}\label{eq:pde}
\left\{
\begin{aligned}
\rho_f\partial_t\bu-\nabla\cdot\bs_f(\bu,p_f)&=\boldsymbol f_f,&&\text{in }\Omega_f,\\
\nabla \cdot \bu&=0, &&\text{in }\Omega_f,\\
\rho_p\partial_t\bxi-\nabla\cdot\bs_p(\be,\phi)&=\boldsymbol f_p,&&\text{in }\Omega_p,\\
\partial_t\be&=\bxi, &&\text{in }\Omega_p,\\
C_0\partial_t\phi+\alpha\nabla\cdot\bxi-\nabla\cdot(K\nabla\phi)&=g_p,
&&\text{in }\Omega_p,\\
\bu_p&=-K\nabla\phi, &&\text{in }\Omega_p.
\end{aligned}
\right.
\end{equation}
Here $\rho_f,\rho_p,\mu_f,\mu_p,C_0>0$, $\lambda_p\ge0$, and the Biot--Willis coefficient $\alpha \in\mathbb R$ is fixed \cite{WangLeeCanic2026SecondOrder}. The time-independent conductivity $K \in L^\infty(\Omega_p)^{d\times d}$ is symmetric and uniformly positive definite. We impose the homogeneous Dirichlet boundary conditions, namely $\bu=0$ on $\Sigma_f$, $\be=\bxi=0$ on $\Sigma_p$, and $\phi=0$ on $\Sigma_p$, and prescribe compatible initial data. We note that these boundary conditions are part of the analytical setting, rather than a claim covering every application boundary condition.

\subsection{Coupling conditions and energy}
On the fluid structure interface $\Gamma$, we consider the following physical coupling conditions, as reported in (cf.~\cite{WangLeeCanic2026SecondOrder})
\begin{subequations}\label{eq:physical}
\begin{align}
\bu\cdot\bn_f&=(\bxi- K \nabla \phi )\cdot\bn_f,\label{eq:flux}\\
\bP\bs_f\bn_f&=-\gamma\bP(\bu-\bxi),\label{eq:bjs}\\
\bn_f\cdot\bs_f\bn_f&=-\phi,\label{eq:normalfluid}\\
\bs_f\bn_f&=-\bs_p\bn_p.\label{eq:balance}
\end{align}
\end{subequations}
Where $\gamma>0$ is the Beavers–Joseph–Saffman  coefficient. In particular, we can derive the mechanical normal condition on the porous side from \eqref{eq:normalfluid} and \eqref{eq:balance}, satisfying
\begin{equation}
\bn_p\cdot\bs_p(\be,\phi)\bn_p+\phi=0.
\label{eq:neumann}
\end{equation}
It is Neumann-type with respect to the mechanics equation and remains coupled to the unknown pore pressure.

We consider the following function spaces
\begin{align*}
\Vf&=\{\bv\in H^1(\Omega_f)^d:\bv|_{\Sigma_f}=0\},
&Q_f&=L^2(\Omega_f),\\
\Vp&=\{\bz\in H^1(\Omega_p)^d:\bz|_{\Sigma_p}=0\},
&Q_p&=\{\psi\in H^1(\Omega_p):\psi|_{\Sigma_p}=0\}.
\end{align*}

For  inner products $(\cdot,\cdot)_i$ on $L^2(\Omega_i)$ ($i=f,p$)
and interface pairing $\langle\cdot,\cdot\rangle_\Gamma$, we introduce the following bilinear forms and norm:
\begin{align*}
a_f(\bu,\bv)&=2\mu_f(\bD(\bu),\bD(\bv))_f,\\
a_e(\be,\bz)&=2\mu_p(\bD(\be),\bD(\bz))_p+
\lambda_p(\nabla\cdot\be,\nabla\cdot\bz)_p,\\
d(\phi,\psi)&=(K\nabla\phi,\nabla\psi)_p,
\qquad \norm{\be}{a_e}^2=a_e(\be,\be).
\end{align*}
The coupled weak formulation of \eqref{eq:pde}, with $\partial_t\be=\bxi$,  is
\begin{subequations}\label{eq:weakphysical}
\begin{align}
&\rho_f(\partial_t\bu,\bv)_f+a_f(\bu,\bv)
-(p_f,\nabla\cdot\bv)_f+(\nabla\cdot\bu,q)_f
+\dual{\phi}{\bv\cdot\bn_f}
+\gamma\dual{\bP(\bu-\bxi)}{\bP\bv}
=(\boldsymbol f_f,\bv)_f,\label{eq:weakfluid}\\
&\rho_p(\partial_t\bxi,\bz)_p+a_e(\be,\bz)
-\alpha(\phi,\nabla\cdot\bz)_p
+\dual{\phi}{\bz\cdot\bn_p}+\gamma\dual{\bP(\bxi-\bu)}{\bP\bz}
=(\boldsymbol f_p,\bz)_p,\label{eq:weakmech}\\
&C_0(\partial_t\phi,\psi)_p+\alpha(\nabla\cdot\bxi,\psi)_p
+d(\phi,\psi)-\dual{\bxi\cdot\bn_p}{\psi}
+\dual{\bu\cdot\bn_p}{\psi}
=(g_p,\psi)_p.\label{eq:weakdarcy}
\end{align}
\end{subequations}
Taking $(\bv,q,\bz,\psi)=(\bu,p_f,\bxi,\phi)$ as the test functions
cancels both the volume Biot coupling terms and the paired normal
pore-pressure and mechanical terms. The normal fluid--Darcy exchange
terms also cancel because the interface normals are oppositely
oriented. Hence, for smooth solutions,
 \DIFaddbegin \begin{equation}
\begin{aligned}% REV: split into two lines (energy written out)
\frac{d}{dt}\frac12\Big(\rho_f\norm{\bu}{f}^2+\rho_p\norm{\bxi}{p}^2
+\norm{\be}{a_e}^2+C_0\norm{\phi}{p}^2\Big)
&+2\mu_f\norm{\bD(\bu)}{f}^2
+\norm{K^{1/2}\nabla\phi}{p}^2
+\gamma\norm{\bP(\bu-\bxi)}{\Gamma}^2\\
&=(\boldsymbol f_f,\bu)_f+(\boldsymbol f_p,\bxi)_p+(g_p,\phi)_p. 
\end{aligned}
\label{eq:physicalenergy}
\end{equation} \DIFaddend 
 % REV-DEL (notation table, (2.6)): the name E_kinetic, which also contained the elastic and storage energies.
Throughout, $\norm{\cdot}{i}$, $i=f,p,\Gamma$ denotes the $L^2$ norm on the corresponding domain  $\Omega_i$ or interface $\Gamma$.

\section{A  Neumann--Robin splitting  scheme }
\label{sec:scheme}
\subsection{ Modified interface conditions }
For   a fixed  $L>0$,   the fluid-side conditions can be  written as 
\begin{subequations}\label{eq:fluidrobin}
\begin{align}
\bn_f\cdot\bs_f\bn_f+L\bu\cdot\bn_f&=L\bu\cdot\bn_f-\phi,\\
\bP\bs_f\bn_f+\gamma\bP\bu&=\gamma\bP\bxi.
\end{align}
\end{subequations}
 For the same $L>0$, the poroelastic-side conditions can be written as 
 \begin{subequations}\label{eq:porousNR}
\begin{align}
\bn_p\cdot\bs_p\bn_p+\phi&=0,\label{eq:NRnormal}\\
K\nabla\phi\cdot\bn_p+L^{-1}\phi-\bxi\cdot\bn_p
&=L^{-1}\phi-\bu\cdot\bn_p,\label{eq:NRdarcy}\\
\bP\bs_p\bn_p+\gamma\bP\bxi&=\gamma\bP\bu.\label{eq:NRtangent}
\end{align}
\end{subequations}
  The  nonzero right-hand sides  %(\eqref{eq:NRdarcy} and \eqref{eq:NRtangent}) 
 will be extrapolated in splitting algorithms . 
In \eqref{eq:NRnormal}, both the stress and pore pressure belong to the current Biot solve, and thus no current fluid unknown appears .

\subsection{Fully discrete algorithm}
  For finite element discretization we consider conforming,  shape-regular, quasi-uniform meshes   of $\Omega_f$ and $\Omega_p$ which are fitted with the interface $\Gamma$. We also assume  mesh size $h\le1$   throughout this paper. We choose conforming finite element  spaces $\Vfh\times Q_{f,h}\subset\Vf\times Q_f$, $\Vph\subset\Vp$, and $Q_{p,h}\subset Q_p$. Displacement and structure velocity use the same space $\Vph$   for discretization. We assume that there exists $\beta_f >0$ independent of $h$ such that 
\begin{align}
	\label{eq:inf-sup}
	\inf_{0\ne{}q_h \in Q_{f,h}} \sup_{0\ne{}\bv_h \in \boldsymbol{V}_{f,h}} \frac{(\nabla \cdot \bv_h, q_h)_f}{\norm{\bv_h}{H^1(\Omega_f)} \norm{q_h}{f} } \ge \beta_f >0 .
\end{align}
There are well-known finite elements satisfying this condition, such as the MINI element \cite{MINI} and the Taylor--Hood element \cite{Taylor-Hood}; their  inf--sup   stability under mild mesh conditions is analyzed in \cite[Chapter~8]{BoffiBrezziFortin2013}. These results are usually stated for velocities vanishing on the whole boundary and mean-zero pressures. Since the functions in $\Vfh$ vanish only on $\Sigma_f$, \eqref{eq:inf-sup} is required for the full pressure space $Q_{f,h}$. The extension to the full pressure space follows from the mean-zero inf--sup condition and a family of discrete velocities whose $H^1$ norms are uniformly bounded and whose fluxes through $\Gamma$ are bounded away from zero. Such velocities can be obtained, for sufficiently small $h$, by boundary-preserving interpolation of a fixed velocity with nonzero interface flux.
% REV11-DEL: ", such as the MINI and Taylor--Hood elements under mild mesh conditions [MINI,Taylor-Hood,BBF]" and
% "; it follows from the mean-zero case by controlling the mean value of the pressure with a fixed discrete velocity
% that has nonzero flux through Gamma." In this paper we will use the family of Taylor-Hood elements.  

 Let $T>0$ be the final time, let $N$ be a positive integer, and set $\Delta t=T/N$.
 For a uniform temporal grid $t_n=n\Delta t$,  $n=0,\ldots,N$,  define the BDF2 approximation and second-order linear extrapolation (EX2) as 
 \begin{equation}
\dt a^{n+1}=\frac{3a^{n+1}-4a^n+a^{n-1}}{2\Delta t},
\qquad a^{*,n+1}=2a^n-a^{n-1},\qquad n\ge1.
\label{eq:operators}
\end{equation}
 For a time-dependent function $w$, we write $w^n=w(t_n)$. In particular, $\boldsymbol f_f^{n+1}=\boldsymbol f_f(t_{n+1})$, $\boldsymbol f_p^{n+1}=\boldsymbol f_p(t_{n+1})$, and $g_p^{n+1}=g_p(t_{n+1})$. The operators in \eqref{eq:operators} are applied in the same way to sampled exact and projected fields.

  Let
 \begin{align}
a_f^R(\bu,\bv)&=a_f(\bu,\bv)+\gamma\dual{\bP\bu}{\bP\bv}
+L\dual{\bu\cdot\bn_f}{\bv\cdot\bn_f},\label{eq:afR}\\
d^R(\phi,\psi)&=d(\phi,\psi)+L^{-1}\dual{\phi}{\psi}.
\label{eq:dR}
\end{align}
Given levels $n$ and $n-1$, perform the following two solves independently.

\paragraph{Fluid solve.}
Find $(\bu_h^{n+1},p_{f,h}^{n+1})\in\Vfh\times Q_{f,h}$ satisfying
  \begin{align}
\rho_f(\dt\bu_h^{n+1},\bv_h)_f+a_f^R(\bu_h^{n+1},\bv_h)
-(p_{f,h}^{n+1},\nabla\cdot\bv_h)_f+(\nabla\cdot\bu_h^{n+1},q_h)_f\notag\\
=(\boldsymbol f_f^{n+1},\bv_h)_f
+\gamma\dual{\bP\bxi_h^{*,n+1}}{\bP\bv_h}
+\dual{L\bu_h^{*,n+1}\cdot\bn_f-\phi_h^{*,n+1}}{\bv_h\cdot\bn_f}
\label{eq:schemefluid}
\end{align}  
 for $(\bv_h, q_h) \in \Vfh \times Q_{f,h}$.

\paragraph{Biot solve.}
Find $(\be_h^{n+1},\bxi_h^{n+1},\phi_h^{n+1})\in\Vph\times\Vph\times Q_{p,h}$ such that
\begin{subequations}\label{eq:schemebiot}
  \begin{align}
&(\dt\be_h^{n+1}, \boldsymbol{\chi}_h)_p = (\bxi_h^{n+1}, \boldsymbol{\chi}_h)_p,\label{eq:schemeinematic}\\
&\rho_p(\dt\bxi_h^{n+1},\bz_h)_p+a_e(\be_h^{n+1},\bz_h)
-\alpha(\phi_h^{n+1},\nabla\cdot\bz_h)_p
+\gamma\dual{\bP\bxi_h^{n+1}}{\bP\bz_h}
+\dual{\phi_h^{n+1}}{\bz_h\cdot\bn_p}\notag\\
&=(\boldsymbol f_p^{n+1},\bz_h)_p
+\gamma\dual{\bP\bu_h^{*,n+1}}{\bP\bz_h},\label{eq:schememech}\\
&C_0(\dt\phi_h^{n+1},\psi_h)_p
+\alpha(\nabla\cdot\bxi_h^{n+1},\psi_h)_p
+d^R(\phi_h^{n+1},\psi_h)
-\dual{\bxi_h^{n+1}\cdot\bn_p}{\psi_h}\notag\\
&=(g_p^{n+1},\psi_h)_p
+\dual{-\bu_h^{*,n+1}\cdot\bn_p+L^{-1}\phi_h^{*,n+1}}{\psi_h}
\label{eq:schemedarcy}
\end{align}  
\end{subequations}
 for $(\boldsymbol{\chi}_h, \bz_h,\psi_h) \in \Vph \times \Vph \times Q_{p,h}$.

The method exchanges data explicitly between the fluid and Biot
subproblems during each time step, while solving each subproblem implicitly. 
  For error analysis the  starting values at $n=0$ and $n=1$   need to be good approximations of the continuous solutions. 
The error conditions will be explicitly given  in Section~\ref{sec:error}.

 %\subsection{The exact difference from the earlier Robin--Robin scheme}
%The earlier mechanical normal condition has the form
%\begin{equation}
%\bn_p\cdot\bs_p^{n+1}\bn_p+\phi^{n+1}
%+\beta\bxi^{n+1}\cdot\bn_p
%=\beta\bxi^{*,n+1}\cdot\bn_p,
%\label{eq:RRnormal}
%\end{equation}
%with $\beta$ representing the robin boundary parameter \cite{WangLeeCanic2026SecondOrder}. Compared with \eqref{eq:schememech}, its weak equation contains two additional terms, namely
%\begin{equation}
%\beta\dual{\bxi_h^{n+1}\cdot\bn_p}{\bz_h\cdot\bn_p}
%\quad\text{on the left},\qquad
%\beta\dual{\bxi_h^{*,n+1}\cdot\bn_p}{\bz_h\cdot\bn_p}
%\quad\text{on the right}.
%\label{eq:deletedpair}
%\end{equation}
%The new scheme removes both terms. In residual form, the difference is
%\begin{equation}
%\beta\dual{(\bxi_h^{n+1}-2\bxi_h^n+\bxi_h^{n-1})\cdot\bn_p}{\bz_h\cdot\bn_p}.
%\label{eq:normaldifference}
%\end{equation}
%For a smooth exact solution this is $O(\Delta t^2)$, so removing it preserves formal second-order consistency. It nevertheless changes the discrete method. Condition~\eqref{eq:NRnormal} is enforced in the natural variational sense, while a finite element stress reconstructed from element gradients need not satisfy a pointwise traction identity.

\begin{proposition}[Subproblem solvability]
\label{prop:solvability}
  Assuming the  discrete inf--sup   condition \eqref{eq:inf-sup} is true ,
each step of \eqref{eq:schemefluid}--\eqref{eq:schemebiot}
has a unique solution for every $\Delta t>0$ , given previous numerical solutions and the right-hand sides .
\end{proposition}

\begin{proof}
  For fixed previous numerical solutions and right-hand sides , both subproblems are   square finite-dimensional linear systems, so it suffices to prove uniqueness.
% REV10-DEL: "both subproblems are linear, so it is enough to show that they have unique solutions."

For the  difference of two fluid solutions,  say $\delta \bu_{h}$,  testing 
the momentum equation   \eqref{eq:schemefluid} with $\delta \bu_{h}$  and using
the incompressibility constraint   yields $\frac{3\rho_f}{2\Delta t}\norm{\delta\bu_h}f^2+a_f^R(\delta \bu_{h}, \delta \bu_{h}) = 0$ . The positive mass, viscous, and Robin contributions   on the boundary imply that $\delta \bu_{h} = 0$ . 
The discrete inf--sup condition  \eqref{eq:inf-sup}  then implies that the pressure difference also vanishes.

For the Biot subproblem, denote the differences  of two solutions  at time
$t_{n+1}$ by $\delta\be_h$, $\delta\bxi_h$, and $\delta\phi_h$.
Since the two solutions have identical histories, subtraction
of their BDF2 kinematic   equation \eqref{eq:schemeinematic}  gives
$\frac{3}{2\Delta t}\delta\be_h=\delta\bxi_h$ .
Test the difference equation of \eqref{eq:schememech}  with $\delta\bxi_h$
and the   difference equation of \eqref{eq:schemedarcy}  with $\delta\phi_h$ , and add them. 
Using the relation $\frac{3}{2\Delta t}\delta\be_h=\delta\bxi_h$  in the elastic term   then yields
% REV-DEL (item 7, Prop. 3.1): the two displayed cancellation identities
%   -\alpha(\delta\phi_h,\nabla\cdot \delta\bxi_h)_{\Omega_p}+\alpha(\nabla\cdot\delta\bxi_h,\delta\phi_h)_{\Omega_p}=0,
%   \langle\delta\phi_h,\delta\bxi_h\cdot\bn_p\rangle_\Gamma-\langle\delta\bxi_h\cdot\bn_p,\delta\phi_h\rangle_\Gamma=0.
 \begin{align*}
0={}&
\frac{3\rho_p}{2\Delta t}\norm{\delta\bxi_h}p^2
+\frac{2\Delta t}{3}\norm{\delta\bxi_h}{a_e}^2
+\frac{3C_0}{2\Delta t}\norm{\delta\phi_h}p^2
+d^R(\delta\phi_h,\delta\phi_h)
+\gamma\norm{\bP\delta\bxi_h}\Gamma^2.
\end{align*}
 Since $\rho_p>0$ and $C_0>0$, this identity implies $\delta\bxi_h=0$ ,  $\delta\phi_h=0$ , 
and therefore $\delta\be_h=0$.
 \end{proof}

\section{Discrete energy stability under a linear CFL condition}
Since the fluid and porous structure exchange energy through the interface and our proposed scheme uses explicit predicted interface values (EX2), that exchange is not perfectly balanced at the discrete level. In this section, we show that the resulting imbalance can be controlled when the time step is sufficiently small relative to the mesh size.

\label{sec:stability}
\subsection{Trace   inequality  and BDF2 energy}
  It is known that there are  constants $C_f,C_p$ independent of $h$ such that
  \begin{equation}
\norm{\bv_h}\Gamma^2\le C_f\norm{\bD(\bv_h)}f^2,
\qquad
\norm{\psi_h}\Gamma^2\le C_p\norm{K^{1/2}\nabla\psi_h}p^2
\label{eq:regulartrace}
\end{equation}  
which are obtained  from continuous trace inequalities, Korn's inequality, and the   Poincar\'{e} inequality (see, e.g., \cite{DiPietro-Ern-book, Brenner-Scott}). 
 For discrete scalar or vector fields on either subdomain,   we have an inverse trace inequality
 \begin{equation}
\norm{w_h}\Gamma^2\le C_Ih^{-1}\norm{w_h}{i}^2,
\qquad i\in\{f,p\},
\label{eq:inversetrace}
\end{equation}
with   $C_I>0$ independent of $h$ .

For any Hilbert inner product, the BDF2 approximation has the following identity
  \begin{align}
(\dt a^{n+1},a^{n+1})
={}&\frac{1}{4\Delta t}\big(\norm{a^{n+1}}{}^2+\norm{2a^{n+1}-a^n}{}^2
-\norm{a^n}{}^2-\norm{2a^n-a^{n-1}}{}^2\big)\notag\\% REV10: split into two lines
&+\frac{1}{4\Delta t}\norm{a^{n+1}-2a^n+a^{n-1}}{}^2.
\label{eq:BDFidentity}
\end{align}  
The identity applies also to $a_e$. Define
  \begin{align}
\E^n:=\frac14\Big[&\rho_f\big(\norm{\bu_h^n}f^2+\norm{2\bu_h^n-\bu_h^{n-1}}f^2\big)
+\rho_p\big(\norm{\bxi_h^n}p^2+\norm{2\bxi_h^n-\bxi_h^{n-1}}p^2\big)\notag\\
&+C_0\big(\norm{\phi_h^n}p^2+\norm{2\phi_h^n-\phi_h^{n-1}}p^2\big)
+\norm{\be_h^n}{a_e}^2+\norm{2\be_h^n-\be_h^{n-1}}{a_e}^2\Big].
\label{eq:discreteenergy}
\end{align}  
Let   $r_a^{n+1}:=a_h^{n+1}-2a_h^n+a_h^{n-1}$ and set
\begin{equation}
\N^{n+1}:=\frac{1}{4\Delta t}\left[
\rho_f\norm{\br_u^{n+1}}f^2+\rho_p\norm{\br_\xi^{n+1}}p^2
+C_0\norm{r_\phi^{n+1}}p^2+\norm{\br_\eta^{n+1}}{a_e}^2\right],
\label{eq:numericaldiss}
\end{equation}  
which corresponds to numerical dissipation.

\subsection{ Discrete stability }
 % REV4: the former Sections 4.2 (The interface residual identity) and 4.3 (Residual bounds and
% stability theorem) are merged into the proof of Theorem 4.1. Only the constants q_Gamma, A, B, c_*
% that enter the statement are defined before it; the moved text is otherwise unchanged.
 Define the constants used in the stability estimate.
%\fourmark{moved from former \S4.3}
 Choose
  \begin{equation}
q_\Gamma
:=\min\left\{\frac{\gamma}{4},\frac{\mu_f}{4C_f}\right\}>0,
\label{eq:qGamma}
\end{equation}  
so that
\[
2q_\Gamma\le\frac{\gamma}{2},
\qquad
2C_fq_\Gamma\le\frac{\mu_f}{2}.
\]
 Define
\begin{equation}
\begin{aligned}
A:=8C_IC_f
\left(
\frac{\gamma^2}{\rho_p}
+\frac{L^2}{\rho_f}
+\frac{1}{C_0}
\right),\qquad
B:=8C_IC_p
\left(
\frac{1}{\rho_f}
+\frac{1}{L^2C_0}
\right),
\end{aligned}
\label{eq:ABconstants}
\end{equation}  
  and choose
\begin{equation}
C_*:=
\min\left\{
\frac{\rho_fq_\Gamma}{16\gamma^2C_I},
\frac{\mu_f}{2A},
\frac{1}{2B}
\right\}.
\label{eq:cstar}
\end{equation}  
  These constants are independent of $h$ and $\Delta t$,
but may depend on the material parameters, Robin boundary coefficient $L$, the geometry,
and mesh regularity.

% REV8-DEL: "For simplicity of presentation we show a discrete stability with zero sources."
\begin{theorem}[Discrete energy stability]
\label{thm:stability}%\eightmark{Remark 4.2 integrated}
Let the sampled sources satisfy $\boldsymbol f_f^j\in\Vf'$, $\boldsymbol f_p^j\in L^2(\Omega_p)^d$, and $g_p^j\in Q_p'$ for $2\le j\le N$, where $\Vf'$ and $Q_p'$ are the dual spaces of $\Vf$ and $Q_p$ equipped with the $H^1$ norms, and the source terms in \eqref{eq:schemefluid} and \eqref{eq:schemebiot} are understood as duality pairings.
If $\Delta t\le C_*h$ for $C_*$ in \eqref{eq:cstar}, then, for $2\le M\le N$, the solution of \eqref{eq:schemefluid}--\eqref{eq:schemebiot} satisfies
\begin{equation}
\begin{aligned}
\E^M+\frac{\Delta t}{2}\sum_{j=2}^{M}\big(\N^j+\Diss_*^j\big)
\le{}&2\E^1+C\Delta t\sum_{j=2}^{M}\Big(\|\boldsymbol f_f^j\|_{\Vf'}^2+\|g_p^j\|_{Q_p'}^2\Big) +\frac{4}{\rho_p}\Big(\Delta t\sum_{j=2}^{M}\norm{\boldsymbol f_p^j}p\Big)^2,
\end{aligned}
\label{eq:stabilitybound}
\end{equation}  
  with $\mathcal{E}^n$, $\mathcal{N}^n$ defined in \eqref{eq:discreteenergy}, \eqref{eq:numericaldiss}, and 
\begin{equation}
\Diss_*^j
:=
\mu_f\norm{\bD(\bu_h^j)}f^2
+\frac12\norm{K^{1/2}\nabla\phi_h^j}p^2
+\frac{q_\Gamma}{2}\norm{\bP\bxi_h^j}\Gamma^2.
\label{eq:Dstar}
\end{equation}  
  Here $C$ depends only on $\mu_f$, the lower bound of the permeability $K$, and the Korn and Poincar\'e constants of $\Vf$ and $Q_p$ but is independent of $h$, $\Delta t$, and $M$. If the sources vanish, then
\begin{equation}
\E^M
+\Delta t\sum_{j=2}^{M}
\left(\frac12\N^j+\Diss_*^j\right)
\le \E^1.
\label{eq:stabilityhom}
\end{equation}
\end{theorem}
% REV8-DEL (theorem statement): "Consider we solve (schemefluid) and (schemebiot) with zero sources."

 \begin{proof}
Steps~1--3 bound the interface residuals and prove \eqref{eq:stabilityhom}; Step~4 treats nonzero sources and proves \eqref{eq:stabilitybound}.
% REV11-DEL: the step-by-step roadmap "The proof consists of 4 steps. In the first step ... In the last step, ..."
% REV8-DEL: "The proof consists of 3 steps. ... In the last step, we complete the proof of (eq:stabilitybound)."

\emph{Step 1: Discrete energy identity.}
%\fourmark{former \S4.2}
 Since   $\dt\be_h^{n+1}-\bxi_h^{n+1}\in\Vph$, taking $\boldsymbol{\chi}_h=\dt\be_h^{n+1}-\bxi_h^{n+1}$ in \eqref{eq:schemeinematic} gives $\bxi_h^{n+1}=\dt\be_h^{n+1}$. We test the fluid equation \eqref{eq:schemefluid} with $(\bv_h, q_h)=(\bu_h^{n+1},p_{f,h}^{n+1})$, the equations \eqref{eq:schememech} and \eqref{eq:schemedarcy} with $\bz_h=\bxi_h^{n+1}$ and $\psi_h = \phi_h^{n+1}$, and add them. By the previous identity, the elastic term is $a_e(\be_h^{n+1},\bxi_h^{n+1})=a_e(\be_h^{n+1},\dt\be_h^{n+1})$.
% REV8-DEL: "the equations (schemebiot) with $\boldsymbol{\chi}_h = \be_h^{n+1}$, $\bz_h^{n+1}=\bxi_h^{n+1}$, $\psi_h = \phi_h^{n+1}$"
The volume coupling terms in \eqref{eq:schemefluid}, \eqref{eq:schemebiot}, and the mechanical normal-interface terms in \eqref{eq:schemebiot} cancel.
The extrapolation in \eqref{eq:operators} satisfies
% REV10-DEL: "Since our scheme use the extrapolation (eq:operators), its difference from the current numerical value is"
\[
a_h^{n+1}-a_h^{*,n+1}
=a_h^{n+1}-2a_h^n+a_h^{n-1}
=r_a^{n+1} \qquad \text{for } a = \boldsymbol{u}, \boldsymbol{\xi}, \phi .
\]  
  Using the identity \eqref{eq:BDFidentity}, and replacing each predictor by $a_h^{*,n+1}=a_h^{n+1}-r_a^{n+1}$, gives
\begin{align}
&\frac{\E^{n+1}-\E^n}{\Delta t}+\N^{n+1}
+2\mu_f\norm{\bD(\bu_h^{n+1})}f^2
+\norm{K^{1/2}\nabla\phi_h^{n+1}}p^2
+\gamma\norm{\bP(\bu_h^{n+1}-\bxi_h^{n+1})}\Gamma^2\notag\\% REV8: split into two lines
&\qquad=\R_\Gamma^{n+1}{}+\mathcal W^{n+1}
\label{eq:exactenergy}
\end{align}  
 where
\begin{align}
\R_\Gamma^{n+1}:={}&
-\gamma\dual{\bP\br_\xi^{n+1}}{\bP\bu_h^{n+1}}
-\gamma\dual{\bP\br_u^{n+1}}{\bP\bxi_h^{n+1}}
-L\dual{\br_u^{n+1}\cdot\bn_f}{\bu_h^{n+1}\cdot\bn_f}\notag\\
&
+\dual{r_\phi^{n+1}}{\bu_h^{n+1}\cdot\bn_f}
+\dual{\br_u^{n+1}\cdot\bn_p}{\phi_h^{n+1}}
-L^{-1}\dual{r_\phi^{n+1}}{\phi_h^{n+1}},
\label{eq:residualsix}
\\
\mathcal W^{n+1}:=&(\boldsymbol f_f^{n+1},\bu_h^{n+1})_f+(\boldsymbol f_p^{n+1},\bxi_h^{n+1})_p+(g_p^{n+1},\phi_h^{n+1})_p. 
\label{eq:sources}
\end{align}
%\reight{Here $\mathcal W^{n+1}:=(\boldsymbol f_f^{n+1},\bu_h^{n+1})_f+(\boldsymbol f_p^{n+1},\bxi_h^{n+1})_p+(g_p^{n+1},\phi_h^{n+1})_p$ is the work of the sources.}

\emph{Step 2: Bounds for the residual terms.}
% REV-DEL (item 7, Sec. 4.3): separate displays of the inverse trace bound for r_xi,
% the fluid trace bound, the Cauchy--Schwarz step, and Young's inequality.
By  the Cauchy--Schwarz inequality, the inverse trace   inequality
\eqref{eq:inversetrace} for $\br_\xi^{n+1}\in\Vph$, the trace estimate
\eqref{eq:regulartrace} for $\bu_h^{n+1}$, and  Young's inequality
  $ab\le\varepsilon a^2+b^2/(4\varepsilon)$ with
$\varepsilon=\rho_p/(32\Delta t)$, we obtain
\begin{align}
\gamma
\left|\dual{\bP\br_\xi^{n+1}}{\bP\bu_h^{n+1}}\right|
&\le
\frac{\rho_p}{32\Delta t}% REV11-DEL: the first of two equal expressions (same terms, coefficients rearranged)
\norm{\br_\xi^{n+1}}p^2
+8\gamma^2C_IC_f\frac{\Delta t}{\rho_ph}
\norm{\bD(\bu_h^{n+1})}f^2.
\label{eq:Rxi}
\end{align}  
 Similarly, applying \eqref{eq:regulartrace}--\eqref{eq:inversetrace}
and Young's inequality yields
  \begin{align}
L\left|
\dual{\br_u^{n+1}\cdot\bn_f}{\bu_h^{n+1}\cdot\bn_f}
\right|
&\le
\frac{\rho_f}{32\Delta t}\norm{\br_u^{n+1}}f^2
+
8L^2C_IC_f\frac{\Delta t}{\rho_fh}
\norm{\bD(\bu_h^{n+1})}f^2,
\label{eq:Runormal}\\
\left|
\dual{r_\phi^{n+1}}{\bu_h^{n+1}\cdot\bn_f}
\right|
&\le
\frac{C_0}{32\Delta t}\norm{r_\phi^{n+1}}p^2
+
8C_IC_f\frac{\Delta t}{C_0h}
\norm{\bD(\bu_h^{n+1})}f^2,
\label{eq:RphiFluid}\\
\left|
\dual{\br_u^{n+1}\cdot\bn_p}{\phi_h^{n+1}}
\right|
&\le
\frac{\rho_f}{32\Delta t}\norm{\br_u^{n+1}}f^2
+
8C_IC_p\frac{\Delta t}{\rho_fh}
\norm{K^{1/2}\nabla\phi_h^{n+1}}p^2,
\label{eq:RuPore}\\
L^{-1}\left|
\dual{r_\phi^{n+1}}{\phi_h^{n+1}}
\right|
&\le
\frac{C_0}{32\Delta t}\norm{r_\phi^{n+1}}p^2
+
8C_IC_p\frac{\Delta t}{L^2C_0h}
\norm{K^{1/2}\nabla\phi_h^{n+1}}p^2.
\label{eq:RphiPore}
\end{align}  
% REV-DEL (item 7, Sec. 4.3): the separate Cauchy--Schwarz display, the
% "insert sqrt(q_Gamma)" explanation, and the two-step rearrangement of the CFL condition.
For the tangential structure-velocity product,   the inverse trace
inequality \eqref{eq:inversetrace} for $\br_u^{n+1}$ and  Young's inequality   give
\begin{equation}
\gamma\left|\dual{\bP\br_u^{n+1}}{\bP\bxi_h^{n+1}}\right|
\le
\frac{q_\Gamma}{2}\norm{\bP\bxi_h^{n+1}}\Gamma^2
+\frac{\gamma^2C_I}{2q_\Gamma h}\norm{\br_u^{n+1}}f^2.
\label{eq:RuTang2}
\end{equation}  
  If
$\Delta t/h\le\rho_fq_\Gamma/(16\gamma^2C_I)$, the last term is at most
$\rho_f\norm{\br_u^{n+1}}f^2/(32\Delta t)$.
Writing $\bP\bxi_h^{n+1}=\bP(\bxi_h^{n+1}-\bu_h^{n+1})+\bP\bu_h^{n+1}$,
using $\|a+b\|^2\le2\|a\|^2+2\|b\|^2$, and applying \eqref{eq:regulartrace}
to $\bu_h^{n+1}$, we obtain from \eqref{eq:qGamma}
\begin{equation}
q_\Gamma\norm{\bP\bxi_h^{n+1}}\Gamma^2
\le
\frac{\gamma}{2}
\norm{\bP(\bu_h^{n+1}-\bxi_h^{n+1})}\Gamma^2
+\frac{\mu_f}{2}\norm{\bD(\bu_h^{n+1})}f^2.
\label{eq:tangentialcontrol}
\end{equation}  
  Notice that this estimate controls the tangential structure velocity through
the slip term and fluid viscous dissipation, without requiring
a bulk $H^1$ estimate for the structure velocity.

  \emph{Step 3: Absorption and summation.}
By the  first restriction in \eqref{eq:cstar} , the condition at the end of Step~2 holds.
 Thus the total coefficients assigned to   $\|\br_u^{n+1}\|_f^2$, $\|\br_{\xi}^{n+1}\|_p^2$, and $\|r_{\phi}^{n+1}\|_p^2$  in
\eqref{eq:Rxi}-- \eqref{eq:RuTang2}  are at most
\[
\frac{3\rho_f}{32\Delta t},\qquad
\frac{\rho_p}{32\Delta t},\qquad
\frac{2C_0}{32\Delta t},
\]
respectively. These contributions are bounded by
$\N^{n+1}/2$.

  Recalling $A, B$ in \eqref{eq:ABconstants} the  remaining restrictions in \eqref{eq:cstar} imply
\[
A\frac{\Delta t}{h}\le\frac{\mu_f}{2},
\qquad
B\frac{\Delta t}{h}\le\frac12.
\]
Combining the residual estimates therefore yields
  \begin{equation}% REV4: numbered (label eq:Rbound) for use in Section 5
\begin{aligned}
|\R_\Gamma^{n+1}|
\le{}&
\frac12\N^{n+1}
+\frac{\mu_f}{2}\norm{\bD(\bu_h^{n+1})}f^2+\frac12\norm{K^{1/2}\nabla\phi_h^{n+1}}p^2
+\frac{q_\Gamma}{2}
\norm{\bP\bxi_h^{n+1}}\Gamma^2.
\end{aligned}
\label{eq:Rbound}
\end{equation}  
Substituting this bound into \eqref{eq:exactenergy}   gives
\begin{equation*}
\begin{aligned}
\frac{\E^{n+1}-\E^n}{\Delta t}
+\frac12\N^{n+1}
&+\frac{3\mu_f}{2}\norm{\bD(\bu_h^{n+1})}f^2
+\frac{1}{2}\norm{K^{1/2}\nabla\phi_h^{n+1}}p^2\\
&+\gamma\norm{\bP(\bu_h^{n+1}-\bxi_h^{n+1})}\Gamma^2
-\frac{q_\Gamma}{2}\norm{\bP\bxi_h^{n+1}}\Gamma^2
\le\mathcal W^{n+1}.
\end{aligned}
\end{equation*}
By \eqref{eq:tangentialcontrol},
\begin{equation*}
\begin{aligned}
&\frac{3\mu_f}{2}\norm{\bD(\bu_h^{n+1})}f^2
+\gamma\norm{\bP(\bu_h^{n+1}-\bxi_h^{n+1})}\Gamma^2
-\frac{q_\Gamma}{2}\norm{\bP\bxi_h^{n+1}}\Gamma^2\\
&\qquad\ge
\mu_f\norm{\bD(\bu_h^{n+1})}f^2
+\frac{q_\Gamma}{2}\norm{\bP\bxi_h^{n+1}}\Gamma^2
+\frac{\gamma}{2}\norm{\bP(\bu_h^{n+1}-\bxi_h^{n+1})}\Gamma^2.
\end{aligned}
\end{equation*}
Inserting this lower bound  and discarding the
remaining nonnegative slip term gives
  \begin{equation}
\begin{aligned}
\frac{\E^{n+1}-\E^n}{\Delta t}
+\frac12\N^{n+1}
&+\mu_f\norm{\bD(\bu_h^{n+1})}f^2
+\frac{1}{2}\norm{K^{1/2}\nabla\phi_h^{n+1}}p^2
+\frac{q_\Gamma}{2}
\norm{\bP\bxi_h^{n+1}}\Gamma^2
\le\mathcal W^{n+1}.
\end{aligned}
\label{eq:stabperstep}
\end{equation}  
  If the sources vanish, then $\mathcal W^{n+1}=0$, and multiplying \eqref{eq:stabperstep}  by $\Delta t$ and summing over $n=1,\ldots,M-1$ proves   \eqref{eq:stabilityhom}.
% REV8-DEL: "Multiplying by Delta t and summing over n=1,...,M-1 proves (eq:stabilitybound)."

  \emph{Step 4: Nonzero sources.}
Since $\bu_h^{n+1}\in\Vf$ and $\phi_h^{n+1}\in Q_p$, Korn's and Poincar\'e's inequalities and Young's inequality give
\begin{equation*}
\begin{aligned}
|(\boldsymbol f_f^{n+1},\bu_h^{n+1})_f|
&\le\frac{\mu_f}{2}\norm{\bD(\bu_h^{n+1})}f^2+C\|\boldsymbol f_f^{n+1}\|_{\Vf'}^2,\\
|(g_p^{n+1},\phi_h^{n+1})_p|
&\le\frac14\norm{K^{1/2}\nabla\phi_h^{n+1}}p^2+C\|g_p^{n+1}\|_{Q_p'}^2.
\end{aligned}
\end{equation*}
These absorptions use one half of the  fluid and pore-pressure   parts of $\Diss_*^{n+1}$, so \eqref{eq:stabperstep} yields
\begin{equation*}
\frac{\E^{n+1}-\E^n}{\Delta t}+\frac12\N^{n+1}+\frac12\Diss_*^{n+1}
\le(\boldsymbol f_p^{n+1},\bxi_h^{n+1})_p
+C\Big(\|\boldsymbol f_f^{n+1}\|_{\Vf'}^2+\|g_p^{n+1}\|_{Q_p'}^2\Big).
\end{equation*}
The mechanical source is tested with the bulk structure velocity, which is not directly controlled by the dissipation, so we bound its time sum through the maximum of the energy. Let $2\le m\le M$ and $\E_{\max}^M:=\max_{1\le j\le M}\E^j$. Since $\rho_p\norm{\bxi_h^j}p^2\le4\E^j$ for $j\ge1$,  Young's inequality   gives
\begin{equation*}
\Delta t\sum_{j=2}^m\big|(\boldsymbol f_p^j,\bxi_h^j)_p\big|
\le\frac{2}{\rho_p^{1/2}}\big(\E_{\max}^M\big)^{1/2}\,\Delta t\sum_{j=2}^M\norm{\boldsymbol f_p^j}p
\le\frac12\E_{\max}^M+\frac{2}{\rho_p}\Big(\Delta t\sum_{j=2}^M\norm{\boldsymbol f_p^j}p\Big)^2.
\end{equation*}
Multiplying the previous per-step inequality by $\Delta t$, summing over $n=1,\ldots,m-1$, and using this bound, we obtain
\begin{equation}
\E^m+\frac{\Delta t}{2}\sum_{j=2}^m\big(\N^j+\Diss_*^j\big)
\le\E^1+\frac12\E_{\max}^M+\mathcal S^M,
\label{eq:stabsum}
\end{equation}
where $\mathcal S^M:=C\Delta t\sum_{j=2}^M\big(\|\boldsymbol f_f^j\|_{\Vf'}^2+\|g_p^j\|_{Q_p'}^2\big)+\frac{2}{\rho_p}\big(\Delta t\sum_{j=2}^M\norm{\boldsymbol f_p^j}p\big)^2$. Inequality \eqref{eq:stabsum} also holds for $m=1$, where the sum on the left is empty. Dropping the nonnegative sum and taking the maximum over $1\le m\le M$ gives $\E_{\max}^M\le\E^1+\frac12\E_{\max}^M+\mathcal S^M$, that is, $\E_{\max}^M\le2\E^1+2\mathcal S^M$. Inserting this into \eqref{eq:stabsum} with $m=M$ gives
$\E^M+\frac{\Delta t}{2}\sum_{j=2}^M(\N^j+\Diss_*^j)\le2\E^1+2\mathcal S^M$, which is \eqref{eq:stabilitybound} after renaming $C$.
\end{proof}

\begin{remark}
The restriction $\Delta t\le C_*h$ is sufficient for the absorption argument above. The optimal admissible time-step constant is not determined here.
% REV10-DEL: "The strict bound C_* in (eq:cstar) is necessary to make the constant C in (eq:stabilitybound)
% independent of M. If Delta t <= C_0 h holds for C_0>0, then we still obtain a discrete stability with a
% discrete Gronwall inequality. Since this argument is rather standard we omit details."
 \end{remark}

 \section{A priori error analysis}
\label{sec:error}

We estimate the difference between the exact solution
$(\bu,p_f,\be,\bxi,\phi)$ and the fully discrete solution
$(\bu_h^n,p_{f,h}^n,\be_h^n,\bxi_h^n,\phi_h^n)$.
The goal is to   prove  an error bound of order $h^k+\Delta t^2$ in the bulk
energy norms under the linear CFL restriction.
  We use the operators in \eqref{eq:operators} and the sampling convention  $w^n=w(t_n)$   of Section~\ref{sec:scheme}.
% REV-DEL (notation table): duplicate definitions of delta_t, w^{*,n+1}, and w^n=w(t_n) (now in Sec. 3.2).

The   analysis is organized as follows. In Section~\ref{subsec:regularity} we introduce the projections of the exact solution and split the errors into projection and discrete parts. In Section~\ref{subsec:errorequations} we derive the equations satisfied by  the discrete errors , whose right-hand sides consist of consistency residuals. These residuals are estimated in Section~\ref{subsec:residualestimates}. Section~\ref{subsec:errortheorem} states the error estimate and proves it by an energy argument, summing in time without a discrete Gronwall inequality.
% REV8-DEL: "The proof is carried out in three main steps. First, we project ... Finally, we bound
% the residuals and sum in time without a discrete Gronwall inequality."
%The only residual requiring a separate argument is the pore-pressure projection load in mechanics. Since the scheme has no normal structure-velocity penalty, its energy estimate does not control
%$\norm{\ehv{\xi}{}\cdot\bn_p}\Gamma$.
%We therefore use the displacement--velocity relation and summation by parts in time to estimate this load through the displacement energy.

\subsection{Regularity assumptions and projection operators}
  \label{subsec:regularity}
 Let $k\ge1$ be an integer. We assume the following regularity
\begin{equation}
\bu,\bxi,\be,\phi\in W^{3,\infty}(0,T;H^{k+1}),
\qquad p_f\in W^{3,\infty}(0,T;H^k),
\qquad \partial_t\be=\bxi,
\label{eq:regularity}
\end{equation}
Here each field is defined on its corresponding subdomain, and vector-valued
spaces are understood.   We define
\begin{equation}
\begin{aligned}
\mathcal M:={}&\norm{\bu}{W^{3,\infty}(0,T;H^{k+1}(\Omega_f))}
+\norm{p_f}{W^{3,\infty}(0,T;H^{k}(\Omega_f))}
+\norm{\be}{W^{3,\infty}(0,T;H^{k+1}(\Omega_p))}\\
&+\norm{\bxi}{W^{3,\infty}(0,T;H^{k+1}(\Omega_p))}
+\norm{\phi}{W^{3,\infty}(0,T;H^{k+1}(\Omega_p))},
\end{aligned}
\label{eq:Mdef}
\end{equation}
where $\norm{w}{W^{3,\infty}(0,T;X)}:=\max_{0\le j\le3}\operatorname*{ess\,sup}_{0<t<T}\norm{\partial_t^jw(t)}{X}$ for a Banach space $X$.
% REV8-DEL: "Let $\mathcal M$ collect the solution norms in (eq:regularity)."
 These assumptions provide the spatial derivatives
needed for an $h^k$ energy estimate and the time derivatives needed for
second-order BDF2 and EX2 consistency. The identity $\be_t=\bxi$ also
supplies the corresponding additional displacement time regularity.
All constants below are independent of mesh size $h$ and time step size $\Delta t$, although they
may depend on the fixed coefficients and the stated regularity bounds. %\revdel
% REV-DEL (notation table): "We explicitly assume K in L^infty(Omega_p)^{d x d}, in addition to the
% stated symmetry and uniform positive definiteness, so that the Darcy bilinear form is continuous
% on H^1(Omega_p)." -- the assumption K in L^infty is now stated in Section 2.

 The projections are chosen to cancel the principal spatial terms when
we subtract the numerical equations. We use a mixed projection for the
fluid and Ritz projections for elasticity and pore pressure.

 At each time, define the Robin--Stokes projection   $(\Piu\bu,\Pip p_f)\in\Vfh\times Q_{f,h}$ by
\begin{align}
a_f^R(\Piu\bu-\bu,\bv_h)-(\Pip p_f-p_f,\nabla\cdot\bv_h)_f&=0,
\qquad (\nabla\cdot(\Piu\bu-\bu),q_h)_f=0.
\label{eq:stokesprojection}
\end{align}  
Here \(a_f^R\) is the fluid bilinear form including the Robin contributions defined in \eqref{eq:afR}.
The  pair $(\Piu\bu,\Pip p_f)$ is defined jointly by \eqref{eq:stokesprojection}, so $\Piu\bu$ and $\Pip p_f$ in general depend on both $\bu$ and $p_f$; we keep the componentwise notation for brevity.
The  divergence condition preserves discrete incompressibility.
 %\revdel
% REV-DEL (item 7, Sec. 5.1): "The first identity cancels the fluid elliptic and pressure-projection
% terms together." -- the cancellations are stated once, after the error equations in Sec. 5.4.

  For projections for the poroelastic subproblem define the  elasticity Ritz map   $\Pieta:\Vp\to\Vph$  and the Darcy Ritz map
  $\Piphi:Q_p\to Q_{p,h}$ by
\begin{equation}
a_e(\Pieta\bz-\bz,\bw_h)=0,\qquad
d^R(\Piphi\psi-\psi,\chi_h)=0
\label{eq:ritzmaps}
\end{equation}  
for all discrete tests. The form $d^R$ is defined in \eqref{eq:dR}.
The projected displacement, structure velocity, and pore pressure are
  $\Pieta\be$, $\Pixi\bxi$, and $\Piphi\phi$, where the velocity projection is
the elasticity Ritz map itself:
\begin{equation}
\Pixi:=\Pieta.
\label{eq:projectedfields}
\end{equation}  
All spaces, coefficients, and maps are time independent.
In particular, the same elasticity map is used for displacement and
velocity, so that
  \[
\partial_t(\Pieta\be)=\Pieta\partial_t\be=\Pieta\bxi=\Pixi\bxi.
\]  
This choice preserves the continuous kinematic relation after spatial
projection and avoids an additional spatial defect in that relation.

 We first state the approximation requirements on the finite element
spaces.   The spaces used in Section~\ref{sec:verification} satisfy them with $k=2$; see the comment after Assumption~\ref{ass:projection}.
% REV10-DEL: "Concrete choices satisfying these requirements are discussed at the end of the section, after the main error argument."

\begin{assumption}[Finite element approximation]\label{ass:projection}% REV8: environment changed from theorem to assumption (cited as Assumption in Lemma 5.2)
In addition to the   stability assumptions above, we  assume the following best-approximation properties, with   constants  independent of $h$:
 % REV10-DEL: "With the inf--sup condition (eq:inf-sup) and the Poincare and Korn inequalities on Omega_p, we can assume the
% following best-approximation properties, with C independent of h:"
 \begin{align}
\inf_{\bv_h\in\Vfh}\norm{\bw-\bv_h}{H^1(\Omega_f)}
+\inf_{q_h\in Q_{f,h}}\norm{r-q_h}f
&\le Ch^k
\left(\norm{\bw}{H^{k+1}(\Omega_f)}+\norm{r}{H^k(\Omega_f)}\right),\notag\\
\inf_{\bz_h\in\Vph}\norm{\bz-\bz_h}{H^1(\Omega_p)}
&\le Ch^k\norm{\bz}{H^{k+1}(\Omega_p)},\notag\\
\inf_{\psi_h\in Q_{p,h}}\norm{\psi-\psi_h}{H^1(\Omega_p)}
&\le Ch^k\norm{\psi}{H^{k+1}(\Omega_p)}.
\label{eq:spaceapproximation}
\end{align}  
The fields in these inequalities belong to the indicated regularity spaces and satisfy their essential boundary conditions.
\end{assumption}
 For the spaces used in Section~\ref{sec:verification}, namely continuous piecewise quadratic elements for the fluid velocity, displacement, structure velocity, and pore pressure and continuous piecewise linear elements for the fluid pressure, \eqref{eq:spaceapproximation} holds with $k=2$ by standard interpolation estimates that preserve the homogeneous boundary conditions on $\Sigma_f$ and $\Sigma_p$ \cite{Brenner-Scott}.

 We split each total error into a spatial projection error and a discrete
error. The former is bounded in Lemma~\ref{lem:projection} below;
the latter will be estimated using the error equations.
The projection error, denoted by the
superscript $I$, is projected minus exact, and the discrete error, denoted by
the superscript $h$, is projected minus discrete. The total errors therefore split as follows:
\begin{subequations}\label{eq:error_decomposition_readable}
\begin{align}
\boldsymbol e_u^n:={}\bu^n-\bu_h^n
&=(\bu^n-\Piu\bu^n)+(\Piu\bu^n-\bu_h^n)
=:-\eIv{u}{n}+\ehv{u}{n},
\label{eq:errsplit_u_new}\\
e_p^n:={}p_f^n-p_{f,h}^n
&=(p_f^n-\Pip p_f^n)+(\Pip p_f^n-p_{f,h}^n)
=:-\eIs{p}{n}+\ehs{p}{n},
\label{eq:errsplit_p_new}\\
\boldsymbol e_\eta^n:={}\be^n-\be_h^n
&=(\be^n-\Pieta\be^n)+(\Pieta\be^n-\be_h^n)
=:-\eIv{\eta}{n}+\ehv{\eta}{n},
\label{eq:errsplit_eta_new}\\
\boldsymbol e_\xi^n:={}\bxi^n-\bxi_h^n
&=(\bxi^n-\Pixi\bxi^n)+(\Pixi\bxi^n-\bxi_h^n)
=:-\eIv{\xi}{n}+\ehv{\xi}{n},
\label{eq:errsplit_xi_new}\\
e_\phi^n:={}\phi^n-\phi_h^n
&=(\phi^n-\Piphi\phi^n)+(\Piphi\phi^n-\phi_h^n)
=:-\eIs{\phi}{n}+\ehs{\phi}{n}.
\label{eq:errsplit_phi_new}
\end{align}
\end{subequations}
The same definitions apply to the time-dependent projected fields before
sampling at $t_n$; for instance, $\eIv{u}{}=\Piu\bu-\bu$ and $\eIs{\phi}{}=\Piphi\phi-\phi$.
% REV2-DEL: the sentence distinguishing the projection error rho_{p_f} from the density rho_p (no longer needed).
The following lemma bounds the projection errors.
% REV10-DEL: "Since Lemma 5.2 controls the projection errors, it remains to estimate the discrete errors through
% their discrete equations."

 \begin{lemma}[Projection estimates sufficient for the energy bound]\label{lem:projection}
Under Assumption~\ref{ass:projection}, the positivity and boundary hypotheses, and \eqref{eq:regularity}, the projections are well defined. For  the projection errors in \eqref{eq:error_decomposition_readable} and  $j=0,1$,
 % REV10-DEL: the repeated definition "Let e_u^I := Pi_u u - u, ..., denote the projection errors."
 \begin{align}
\norm{\partial_t^j\eIv{u}{}}{H^1(\Omega_f)}
 +\norm{\partial_t^j\eIs{p}{}}f
 +\norm{\partial_t^j\eIv{\eta}{}}{H^1(\Omega_p)}
 +\norm{\partial_t^j\eIv{\xi}{}}{H^1(\Omega_p)}
 +\norm{\partial_t^j\eIs{\phi}{}}{H^1(\Omega_p)}
 &\le Ch^k\mathcal M,
\label{eq:approximation}\\
\norm{\partial_t^j\eIv{u}{}}\Gamma
 +\norm{\partial_t^j\eIv{\eta}{}}\Gamma
 +\norm{\partial_t^j\eIv{\xi}{}}\Gamma
 +\norm{\partial_t^j\eIs{\phi}{}}\Gamma
 &\le Ch^k\mathcal M.
\label{eq:projectiontrace}
\end{align}  
  Moreover, for $0\le j\le3$ and almost every $t\in(0,T)$,
\begin{equation}
\begin{aligned}
\norm{\partial_t^j\Piu\bu(t)}{H^1(\Omega_f)}
&\le C\big(\norm{\partial_t^j\bu(t)}{H^1(\Omega_f)}+\norm{\partial_t^jp_f(t)}f\big)\le C\mathcal M,\\
\norm{\partial_t^j\Pieta\be(t)}{H^1(\Omega_p)}+\norm{\partial_t^j\Pixi\bxi(t)}{H^1(\Omega_p)}
&\le C\mathcal M,\qquad
\norm{\partial_t^j\Piphi\phi(t)}{H^1(\Omega_p)}\le C\mathcal M.
\end{aligned}
\label{eq:projtimereg}
\end{equation}
 Both Ritz maps are $H^1$ stable.
 % REV10-DEL: "The four projected fields ... have uniform W^{3,infty}(0,T;H^1) bounds." All constants are independent of $h,\Delta t$; $\mathcal M$ \reight{is defined in \eqref{eq:Mdef}}.
\end{lemma}
\begin{proof}
\emph{Spatial estimates.}
The forms $a_f^R$, $a_e$, and $d^R$ are continuous and coercive in their respective $H^1$ spaces by the continuous trace, Korn, and Poincar\'{e} inequalities   under the stated coefficient and boundary assumptions . The discrete mixed stability estimate for \eqref{eq:stokesprojection} therefore gives
  \[
\norm{\Piu\bu-\bu}{H^1(\Omega_f)}+\norm{\Pip p_f-p_f}f
\le C\left(
\inf_{\bv_h\in\Vfh}\norm{\bu-\bv_h}{H^1(\Omega_f)}
+\inf_{q_h\in Q_{f,h}}\norm{p_f-q_h}f
\right).
\]  
C\'ea's lemma gives the analogous best-approximation estimates for   $\Pieta$ and $\Piphi$ ; see \cite{BoffiBrezziFortin2013}. Combining these estimates with \eqref{eq:spaceapproximation} proves \eqref{eq:approximation} for $j=0$. 

\emph{Time derivatives and interface estimates.}
Since the maps are linear , bounded,  and time independent, they commute with time differentiation , because bounded linear maps commute with Bochner integrals. For the Robin--Stokes projection, $(\partial_t^j\Piu\bu,\partial_t^j\Pip p_f)$ is the projection of the pair $(\partial_t^j\bu,\partial_t^jp_f)$ . Applying the same estimates to the first time derivatives proves the case $j=1$. The continuous trace theorem then gives \eqref{eq:projectiontrace}.

\emph{Stability of the projected fields.}
Coercivity also gives   $\norm{\Pieta\bz}{H^1}\le C\norm{\bz}{H^1}$ and $\norm{\Piphi\psi}{H^1}\le C\norm{\psi}{H^1}$ . Mixed stability yields   $\norm{\Piu\bu}{H^1}+\norm{\Pip p_f}f\le C(\norm{\bu}{H^1}+\norm{p_f}f)$ . Apply these stability bounds to time derivatives through order three to obtain   \eqref{eq:projtimereg}. The bound for $\Piu\bu$ in \eqref{eq:projtimereg} involves the time derivatives of $p_f$; this is why $\mathcal M$ contains the $W^{3,\infty}$ norm of $p_f$.
% REV11-DEL: "Since Pi_u u is the velocity component of the projection of the pair (u,p_f), its time derivatives are
% controlled by those of both u and p_f;" (the joint dependence is explained after (eq:stokesprojection)) In particular, the $W^{3,\infty}(L^2)$ and $W^{2,\infty}(H^1)$ bounds used below follow. No adjoint regularity is used.
 \end{proof}

 % REV8: Sections 5.2--5.6 of v7 are reorganized into: error decomposition (5.2), error equations with
% residuals (5.3), estimates of the residuals (5.4), and the error theorem stated before its proof (5.5).
% Moved text is verbatim unless shown in teal; deletions are recorded in "% REV8-DEL" comments.

  \subsection{Error equations with residuals}
\label{subsec:errorequations}
In this subsection we derive the equations satisfied by the discrete errors. Their right-hand sides consist of consistency defects, which we define first, and of residual functionals built from these defects and the projection errors.

 \paragraph{ Consistency defects .}%\eightmark{from former \S5.3}
Although the exact fields satisfy $\bxi=\partial_t\be$, replacing the
time derivative by BDF2 introduces a mismatch between the projected
velocity and the discrete derivative of projected displacement.
Define
  \begin{equation}
\boldsymbol\kappa^{n+1}
:=\Pixi\bxi^{n+1}-\dt \Pieta\be^{n+1}.
\label{eq:kappa_def}
\end{equation}  
 % REV8-DEL: the subequations wrapper (eq:kinematicerror) and "Since $\bxi_h^{n+1}=\dt\be_h^{n+1}$, subtraction
% gives" with (eq:error_kinematic_relation); the relation is now part of Lemma lem:errorequations.
 The  remaining consistency errors arise from approximating time derivatives
and predicting interface values. Define
  \begin{align*}
\boldsymbol m_u^{n+1}&=\dt \Piu\bu^{n+1}-\partial_t\bu(t_{n+1}),&
\boldsymbol m_\xi^{n+1}&=\dt \Pixi\bxi^{n+1}-\partial_t\bxi(t_{n+1}),
\\
m_\phi^{n+1}&=\dt \Piphi\phi^{n+1}-\partial_t\phi(t_{n+1}),&
\boldsymbol z_u^{n+1}&=\Piu\bu^{*,n+1}-\bu(t_{n+1}),
\\
\boldsymbol z_\xi^{n+1}&=\Pixi\bxi^{*,n+1}-\bxi(t_{n+1}), &
z_\phi^{n+1}&=\Piphi\phi^{*,n+1}-\phi(t_{n+1}).
\end{align*} 

 %\begin{align*}
%\reight{\boldsymbol m_u^{n+1}}&=\dt \Piu\bu^{n+1}-\partial_t\bu(t_{n+1}),\qquad
%\reight{\boldsymbol m_\xi^{n+1}}=\dt \Pixi\bxi^{n+1}-\partial_t\bxi(t_{n+1}),\\
%m_\phi^{n+1}&=\dt \Piphi\phi^{n+1}-\partial_t\phi(t_{n+1}),\qquad
%\boldsymbol z_u^{n+1}=\Piu\bu^{*,n+1}-\bu(t_{n+1}),\\
%\boldsymbol z_\xi^{n+1}&=\Pixi\bxi^{*,n+1}-\bxi(t_{n+1}),\qquad
%z_\phi^{n+1}=\Piphi\phi^{*,n+1}-\phi(t_{n+1}).
%\end{align*}
% REV8-DEL: "Set $\chi=h^k+\Delta t^2$."

  \paragraph{Residual functionals.}%\eightmark{from former \S5.4}
The error  equations have the same principal operators as the numerical
scheme. Their additional right-hand sides record the defects   defined
 above. We group these terms by the equation in which they occur:
  $F_1$  for the fluid,   $F_2$  for elastic mechanics, and   $F_4$  for pore pressure.
The pore-pressure projection load in mechanics is denoted separately by
  $F_3$ , since it requires summation by parts in time.
For compactness, all defects on the right-hand sides below are evaluated
at $t_{n+1}$ unless another time is displayed.   For $\bv_h\in\Vfh$, $\bz_h\in\Vph$, $\bz\in\Vp$, and $\psi_h\in Q_{p,h}$, the residuals are
% REV10-DEL: "Explicit forms of the residuals are:" and the test-function memberships inside the display 
%
\begin{align}
F_1^{n+1}(\bv_h)={}&\rho_f(\boldsymbol m_u,\bv_h)_f
-\gamma\dual{\bP\boldsymbol z_\xi}{\bP\bv_h}
-L\dual{\boldsymbol z_u\cdot\bn_f}{\bv_h\cdot\bn_f}
+\dual{z_\phi}{\bv_h\cdot\bn_f},
\label{eq:F1}
\\
F_2^{n+1}(\bz_h)={}&\rho_p(\boldsymbol m_\xi,\bz_h)_p
+\gamma\dual{\bP(\eIv{\xi}{}-\boldsymbol z_u)}{\bP\bz_h},
\label{eq:F2}
\\
F_3(t;\bz)={}&-\alpha(\eIs{\phi}{}(t),\nabla\cdot\bz)_p
+\dual{\eIs{\phi}{}(t)}{\bz\cdot\bn_p}
\label{eq:F3}
\\
={}&(1-\alpha)(\eIs{\phi}{}(t),\nabla\cdot\bz)_p+(\nabla\eIs{\phi}{}(t),\bz)_p,
\label{eq:F3alt}
\\
F_4^{n+1}(\psi_h)={}&C_0(m_\phi,\psi_h)_p
+\alpha(\nabla\cdot\eIv{\xi}{},\psi_h)_p
-\dual{\eIv{\xi}{}\cdot\bn_p}{\psi_h}
+\dual{\boldsymbol z_u\cdot\bn_p-L^{-1}z_\phi}{\psi_h}.
\label{eq:F4}
\end{align} 
 % REV-DEL (items 6-7, Sec. 5.4): "The residual J tests only the bulk velocity and its tangential trace.
% In contrast, ell contains a normal traction contribution."
%{\color{red}The definition \eqref{eq:elldef} extends to all $\bz\in\Vp$. Since
%$\bz|_{\Sigma_p}=0$ and $\eIs{\phi}{}(t)\in H^1(\Omega_p)$, the divergence theorem gives
%$(\eIs{\phi}{}(t),\nabla\cdot\bz)_p=-(\nabla\eIs{\phi}{}(t),\bz)_p+\dual{\eIs{\phi}{}(t)}{\bz\cdot\bn_p}$,
%and hence the equivalent form without interface terms}
%\begin{equation}\color{red}
%\reight{F_3}(t;\bz)=(1-\alpha)(\eIs{\phi}{}(t),\nabla\cdot\bz)_p+(\nabla\eIs{\phi}{}(t),\bz)_p,
%\qquad \bz\in\Vp.
%\label{eq:F3alt}
%\end{equation}
%\begin{align}
%\reight{F_4^{n+1}}(\psi_h)={}&C_0(m_\phi,\psi_h)_p
%+\alpha(\nabla\cdot\eIv{\xi}{},\psi_h)_p
%-\dual{\eIv{\xi}{}\cdot\bn_p}{\psi_h}
%+\dual{\boldsymbol z_u\cdot\bn_p-L^{-1}z_\phi}{\psi_h}.
%\label{eq:Gdef}
%\end{align}
%The normal structure trace here belongs to the known projection error $\eIv{\xi}{}$, which is controlled by \eqref{eq:projectiontrace}.

 %\paragraph{Error equations.}
 \begin{lemma}[Error equations]\label{lem:errorequations}
For $1\le n\le N-1$, the discrete errors satisfy the kinematic relation
\begin{equation}
\ehv{\xi}{n+1}
=\dt \ehv{\eta}{n+1}+\boldsymbol\kappa^{n+1},
\label{eq:error_kinematic_relation}
\end{equation}  
  the fluid error equation
\begin{align}
&\rho_f(\dt \ehv{u}{n+1},\bv_h)_f+a_f^R(\ehv{u}{n+1},\bv_h)
-(\ehs{p}{n+1},\nabla\cdot\bv_h)_f
+(\nabla\cdot\ehv{u}{n+1},q_h)_f\notag\\
&=\gamma\dual{\bP\ehv{\xi}{*,n+1}}{\bP\bv_h}
+\dual{L\ehv{u}{*,n+1}\cdot\bn_f-\ehs{\phi}{*,n+1}}{\bv_h\cdot\bn_f}+F_1^{n+1}(\bv_h),
\label{eq:errorfluid}
\end{align}  
  for all $(\bv_h,q_h)\in\Vfh\times Q_{f,h}$, and the Biot error equations
 \begin{subequations}\label{eq:errorbiot}
  \begin{align}
&\rho_p(\dt \ehv{\xi}{n+1},\bz_h)_p+a_e(\ehv{\eta}{n+1},\bz_h)
-\alpha(\ehs{\phi}{n+1},\nabla\cdot\bz_h)_p+\gamma\dual{\bP\ehv{\xi}{n+1}}{\bP\bz_h}
+\dual{\ehs{\phi}{n+1}}{\bz_h\cdot\bn_p}\notag\\
&=\gamma\dual{\bP\ehv{u}{*,n+1}}{\bP\bz_h}
+F_2^{n+1}(\bz_h)+F_3(t_{n+1};\bz_h),\label{eq:errormechanics}\\
&C_0(\dt \ehs{\phi}{n+1},\psi_h)_p
+\alpha(\nabla\cdot\ehv{\xi}{n+1},\psi_h)_p
+d^R(\ehs{\phi}{n+1},\psi_h)
-\dual{\ehv{\xi}{n+1}\cdot\bn_p}{\psi_h}\notag\\
&=\dual{-\ehv{u}{*,n+1}\cdot\bn_p+L^{-1}\ehs{\phi}{*,n+1}}{\psi_h}
+F_4^{n+1}(\psi_h)
\label{eq:errordarcy}
\end{align}  
\end{subequations}
  for all $(\bz_h,\psi_h)\in\Vph\times Q_{p,h}$.
\end{lemma}
% REV8-DEL: "Subtracting the numerical equations from the projected continuous equations gives" (now the proof below).
\begin{proof}
%\emph{Kinematic relation.} 
Since $\dt\be_h^{n+1}-\bxi_h^{n+1}\in\Vph$, the kinematic equation \eqref{eq:schemeinematic} gives $\bxi_h^{n+1}=\dt\be_h^{n+1}$. Subtracting this identity from \eqref{eq:kappa_def} gives \eqref{eq:error_kinematic_relation}.

%\emph{Fluid equation.} 
For the fluid equation we evaluate \eqref{eq:weakfluid} at $t_{n+1}$ and test it with $(\bv_h,q_h)$. Adding $L\dual{\bu\cdot\bn_f}{\bv_h\cdot\bn_f}$ to both sides and moving the terms containing $\bxi$  and  $\phi$ to the right,  the   left-hand side contains $a_f^R(\bu,\bv_h)$ from \eqref{eq:afR}. By \eqref{eq:stokesprojection}, the elliptic, pressure, and divergence terms are unchanged when $(\bu,p_f)$ is replaced by $(\Piu\bu,\Pip p_f)$. Inserting $\partial_t\bu(t_{n+1})=\dt\Piu\bu^{n+1}-\boldsymbol m_u^{n+1}$ and, in the interface terms, $\bu(t_{n+1})=\Piu\bu^{*,n+1}-\boldsymbol z_u^{n+1}$, $\bxi(t_{n+1})=\Pixi\bxi^{*,n+1}-\boldsymbol z_\xi^{n+1}$, and $\phi(t_{n+1})=\Piphi\phi^{*,n+1}-z_\phi^{n+1}$, we obtain
\begin{equation*}
\begin{aligned}
&\rho_f(\dt\Piu\bu^{n+1},\bv_h)_f+a_f^R(\Piu\bu^{n+1},\bv_h)
-(\Pip p_f^{n+1},\nabla\cdot\bv_h)_f+(\nabla\cdot\Piu\bu^{n+1},q_h)_f\\
&=(\boldsymbol f_f^{n+1},\bv_h)_f+\gamma\dual{\bP\Pixi\bxi^{*,n+1}}{\bP\bv_h}
+\dual{L\Piu\bu^{*,n+1}\cdot\bn_f-\Piphi\phi^{*,n+1}}{\bv_h\cdot\bn_f}+F_1^{n+1}(\bv_h).
\end{aligned}
\end{equation*}
Subtracting \eqref{eq:schemefluid} and using the linearity of the extrapolation, for instance $\Piu\bu^{*,n+1}-\bu_h^{*,n+1}=\ehv{u}{*,n+1}$, gives \eqref{eq:errorfluid}.

%\emph{Biot equations.} 
For the Biot equation, by \eqref{eq:ritzmaps}, $a_e(\be,\bz_h)=a_e(\Pieta\be,\bz_h)$ and $d^R(\phi,\psi_h)=d^R(\Piphi\phi,\psi_h)$. We evaluate \eqref{eq:weakmech} at $t_{n+1}$, test it with $\bz_h$, and insert $\partial_t\bxi(t_{n+1})=\dt\Pixi\bxi^{n+1}-\boldsymbol m_\xi^{n+1}$, $\bxi^{n+1}=\Pixi\bxi^{n+1}-\eIv{\xi}{n+1}$, $\phi^{n+1}=\Piphi\phi^{n+1}-\eIs{\phi}{n+1}$, and, in the tangential coupling term, $\bu(t_{n+1})=\Piu\bu^{*,n+1}-\boldsymbol z_u^{n+1}$. The result is \eqref{eq:schememech} written for the projected fields, with the additional terms $F_2^{n+1}(\bz_h)+F_3(t_{n+1};\bz_h)$ on the right-hand side; subtracting \eqref{eq:schememech} gives \eqref{eq:errormechanics}. Similarly, we evaluate \eqref{eq:weakdarcy} at $t_{n+1}$, test it with $\psi_h$, and add $L^{-1}\dual{\phi}{\psi_h}$ to both sides, so that $d^R$ from \eqref{eq:dR} appears on the left. Inserting $\partial_t\phi(t_{n+1})=\dt\Piphi\phi^{n+1}-m_\phi^{n+1}$ and $\bxi^{n+1}=\Pixi\bxi^{n+1}-\eIv{\xi}{n+1}$ on the left, and $\bu(t_{n+1})=\Piu\bu^{*,n+1}-\boldsymbol z_u^{n+1}$ and $\phi(t_{n+1})=\Piphi\phi^{*,n+1}-z_\phi^{n+1}$ on the right, gives \eqref{eq:schemedarcy} written for the projected fields, with the additional term $F_4^{n+1}(\psi_h)$. Subtracting \eqref{eq:schemedarcy} gives \eqref{eq:errordarcy}.
\end{proof}
%No term $\dual{\ehv{\xi}{n+1}\cdot\bn_p}{\bz_h\cdot\bn_p}$ or its extrapolated counterpart occurs in these equations. The remaining normal pressure-projection contribution is explicitly retained in \eqref{eq:elldef}.

% REV8-DEL: "The Stokes projection cancels the fluid elliptic and pressure-projection pair, the elasticity
% projection cancels the displacement elliptic error, and the Darcy projection cancels the scalar d^R error."
% (now shown in the proof of Lemma lem:errorequations)

\subsection{Estimates of the residuals}
\label{subsec:residualestimates}
We now estimate the defects and residual functionals appearing in Lemma~\ref{lem:errorequations}.
% REV8-DEL: "Throughout, we set $\chi=h^k+\Delta t^2$." (chi is no longer used)

\begin{lemma}[Consistency defects]\label{lem:defects}%\eightmark{from former \S5.3}
For $1\le n\le N-1$,
\begin{equation}
\norm{\boldsymbol\kappa^{n+1}}{H^1(\Omega_p)}\le C\Delta t^2\mathcal M.
\label{eq:kappabound}
\end{equation}
Moreover, the derivative defects are measured in bulk
$L^2$ norms, while the extrapolation defects are measured on the interface:
\begin{align}
\norm{\boldsymbol m_u^{n+1}}f+
\norm{\boldsymbol m_\xi^{n+1}}p+
\norm{m_\phi^{n+1}}p&\le C(h^k+\Delta t^2)\mathcal M,
\label{eq:consistencybounds}
\\
\norm{\boldsymbol z_u^{n+1}}\Gamma+
\norm{\boldsymbol z_\xi^{n+1}}\Gamma+
\norm{z_\phi^{n+1}}\Gamma&\le C(h^k+\Delta t^2)\mathcal M.
\label{eq:defectbounds}
\end{align}
\end{lemma}
\begin{proof}
% REV8-DEL: previous short proof ("Taylor's formula gives ...", "The first term is bounded by C dt^2 M in
% L^2 using the third time derivative of Pi_u u ...", "The same argument applies to m_p and m_phi",
% "The bounds for z_xi and z_phi follow identically. Thus both sources of consistency error have the
% common size chi."); replaced by the detailed proof below.
%\emph{Taylor estimates.} 
We first show two consistency error estimates using Taylor expansions.
Let $X$ be a Banach space, $w\in W^{3,\infty}(0,T;X)$, and $I_n=(t_{n-1},t_{n+1})$. Taylor's formula with integral remainder at $t_{n+1}$ gives, for $s\in\{\Delta t,2\Delta t\}$,
\[
\begin{aligned}
w(t_{n+1}-s)&=w(t_{n+1})-s\,\partial_tw(t_{n+1})+\frac{s^2}{2}\,\partial_t^2w(t_{n+1})-R(s),\\
R(s)&=\frac12\int_{t_{n+1}-s}^{t_{n+1}}(\tau-t_{n+1}+s)^2\,\partial_t^3w(\tau)\,d\tau,
\end{aligned}
\]
with $\norm{R(s)}X\le\frac{s^3}{6}\norm{\partial_t^3w}{L^\infty(I_n;X)}$. Hence $3w^{n+1}-4w^n+w^{n-1}=2\Delta t\,\partial_tw(t_{n+1})+4R(\Delta t)-R(2\Delta t)$. Moreover, $w^{n+1}-2w^n+w^{n-1}=\int_{t_{n-1}}^{t_{n+1}}(\Delta t-|\tau-t_n|)\,\partial_t^2w(\tau)\,d\tau$, and the kernel has integral $\Delta t^2$. Therefore
\begin{equation}
\norm{\dt w^{n+1}-\partial_tw(t_{n+1})}X\le\Delta t^2\norm{\partial_t^3w}{L^\infty(I_n;X)},
\qquad
\norm{w^{*,n+1}-w(t_{n+1})}X\le\Delta t^2\norm{\partial_t^2w}{L^\infty(I_n;X)}.
\label{eq:taylorBDF}
\end{equation}

% REV10-DEL: the paragraphs on time derivatives of the projections with (eq:ritztime) and (eq:stokesprojtime);
% these bounds are now stated once, in (eq:projtimereg) of Lemma lem:projection, with the p_f-dependence of Pi_u u.

We now show the main estimates. 
%\emph{Kinematic defect.} 
Since $\Pixi=\Pieta$ is linear and $\bxi=\partial_t\be$, the common elasticity map gives the exact identity
\[
\boldsymbol\kappa^{n+1}
=\Pieta\big(\partial_t\be(t_{n+1})-\dt \be^{n+1}\big).
\]
%\revdel%
% REV-DEL (item 7, Sec. 5.3): "Because the same map R_e is used for both fields, kappa is only a
% time-discretization defect." -- the advantage of the common map is stated once in Sec. 5.1.
The first estimate in \eqref{eq:taylorBDF} with $X=H^1(\Omega_p)^d$ and $w=\be$ gives
\[
\norm{\be_t(t_{n+1})-\dt \be^{n+1}}{H^1(\Omega_p)}
\le C\Delta t^2
\norm{\partial_t^3\be}{L^\infty(t_{n-1},t_{n+1};H^1(\Omega_p))}.
\]
 The   $H^1$ stability of $\Pieta$ therefore yields
\eqref{eq:kappabound}.

%\reight{\emph{Derivative defects.}} 
To show \eqref{eq:consistencybounds}, we split $\boldsymbol m_u^{n+1}$ into a BDF2 defect and a projection error:
\[
\boldsymbol m_u^{n+1}
=\dt \Piu\bu^{n+1}-\partial_t\Piu\bu(t_{n+1})
+\partial_t\eIv{u}{}(t_{n+1}).
% REV11-DEL: \underbrace annotations "BDF2 defect" and "projection error"
\]
By the first estimate in \eqref{eq:taylorBDF} with $X=L^2(\Omega_f)^d$ and $w=\Piu\bu$, and by \eqref{eq:projtimereg} with $j=3$,
\[
\begin{aligned}
\norm{\dt \Piu\bu^{n+1}-\partial_t\Piu\bu(t_{n+1})}f
&\le\Delta t^2\norm{\partial_t^3\Piu\bu}{L^\infty(I_n;H^1(\Omega_f))}\\
&\le C\Delta t^2\Big(\norm{\partial_t^3\bu}{L^\infty(I_n;H^1(\Omega_f))}+\norm{\partial_t^3p_f}{L^\infty(I_n;L^2(\Omega_f))}\Big)
\le C\Delta t^2\mathcal M.
\end{aligned}
\]
The projection error satisfies $\norm{\partial_t\eIv{u}{}(t_{n+1})}f\le Ch^k\mathcal M$ by \eqref{eq:approximation} with $j=1$. For $\boldsymbol m_\xi$ and $m_\phi$ we use the same splitting with $\Pixi\bxi$ and $\Piphi\phi$ in place of $\Piu\bu$. By \eqref{eq:taylorBDF} and \eqref{eq:projtimereg} with $j=3$, their BDF2 defects are bounded by $C\Delta t^2\mathcal M$ in $L^2(\Omega_p)$, and their projection errors $\partial_t\eIv{\xi}{}(t_{n+1})$ and $\partial_t\eIs{\phi}{}(t_{n+1})$ are bounded by $Ch^k\mathcal M$ by \eqref{eq:approximation} with $j=1$. This proves \eqref{eq:consistencybounds}.

%\reight{\emph{Extrapolation defects.}}
For \eqref{eq:defectbounds} we split $\boldsymbol z_u^{n+1}$ into an EX2 defect and a projection error:
\[
\boldsymbol z_u^{n+1}
=\Piu\bu^{*,n+1}-\Piu\bu^{n+1}
+\eIv{u}{n+1}.
% REV11-DEL: \underbrace annotations "EX2 defect" and "projection error"; "consider the splitting"
\]
By the second estimate in \eqref{eq:taylorBDF} with $X=H^1(\Omega_f)^d$ and $w=\Piu\bu$, and by \eqref{eq:projtimereg} with $j=2$, the EX2 defect is bounded in $H^1(\Omega_f)$ by $C\Delta t^2\big(\norm{\partial_t^2\bu}{L^\infty(I_n;H^1(\Omega_f))}+\norm{\partial_t^2p_f}{L^\infty(I_n;L^2(\Omega_f))}\big)\le C\Delta t^2\mathcal M$; the continuous trace theorem gives the same bound in $L^2(\Gamma)$. The projection error satisfies $\norm{\eIv{u}{n+1}}\Gamma\le Ch^k\mathcal M$ by \eqref{eq:projectiontrace}. The bounds for $\boldsymbol z_\xi$ and $z_\phi$ follow in the same way, with the bounds for $\Pixi\bxi$ and $\Piphi\phi$ in \eqref{eq:projtimereg} and the trace theorem on $\Omega_p$. Adding these estimates gives \eqref{eq:defectbounds}.
\end{proof}

To measure the  pressure-projection   load $F_3$, we give $\Vp$ the equivalent elasticity norm
$\norm{\bz}{a_e}=a_e(\bz,\bz)^{1/2}$. Its dual norm is
\[
\|F_3(t)\|_{a_e'}
=\sup_{0\ne\bz\in\Vp}
\frac{|F_3(t;\bz)|}{\norm{\bz}{a_e}}.
\]

  \begin{lemma}[Residual bounds]\label{lem:residualbounds}
%\eightmark{from former \S\S5.4--5.5}
For $1\le n\le N-1$ and all $\bv_h\in\Vfh$, $\bz_h\in\Vph$, and $\psi_h\in Q_{p,h}$,
\begin{align}
|F_1^{n+1}(\bv_h)|&\le C(h^k+\Delta t^2)\mathcal M\norm{\bv_h}{H^1(\Omega_f)},\notag\\
|F_2^{n+1}(\bz_h)|&\le C(h^k+\Delta t^2)\mathcal M\big(\norm{\bz_h}p+\norm{\bP\bz_h}\Gamma\big),\notag\\
|F_4^{n+1}(\psi_h)|&\le C(h^k+\Delta t^2)\mathcal M\norm{\psi_h}{H^1(\Omega_p)}.
\label{eq:FGJbounds}
\end{align}
Moreover, the pressure-projection load satisfies
\begin{equation}
\sup_{0\le t\le T}\|F_3(t)\|_{a_e'}+
\int_0^T\|\partial_tF_3(t)\|_{a_e'}\,dt
\le C_T h^k\mathcal M.
\label{eq:ellbound}
\end{equation}
\end{lemma}
\begin{proof}
 Applying Cauchy--Schwarz and the continuous trace theorem to the
individual terms   of \eqref{eq:F1}, \eqref{eq:F2}, and \eqref{eq:F4} gives \eqref{eq:FGJbounds}.
% REV8-DEL: "Applying Cauchy--Schwarz and the continuous trace theorem to the individual terms gives"
 Here we used \eqref{eq:approximation}--\eqref{eq:projectiontrace}
and \eqref{eq:consistencybounds} --\eqref{eq:defectbounds} . For example,
  \[
\begin{aligned}
|F_2^{n+1}(\bz_h)|
&\le \rho_p\norm{\boldsymbol m_\xi}p\norm{\bz_h}p
+\gamma\big(\norm{\bP\eIv{\xi}{}}\Gamma
+\norm{\bP\boldsymbol z_u}\Gamma\big)
\norm{\bP\bz_h}\Gamma.
\end{aligned}
\]  
  Similarly, since $\norm{\bv_h}\Gamma\le C\norm{\bv_h}{H^1(\Omega_f)}$ and $\norm{\psi_h}\Gamma\le C\norm{\psi_h}{H^1(\Omega_p)}$,
\begin{equation*}
\begin{aligned}
|F_1^{n+1}(\bv_h)|&\le\rho_f\norm{\boldsymbol m_u}f\norm{\bv_h}f
+\big(\gamma\norm{\bP\boldsymbol z_\xi}\Gamma+L\norm{\boldsymbol z_u}\Gamma+\norm{z_\phi}\Gamma\big)\norm{\bv_h}\Gamma,\\
|F_4^{n+1}(\psi_h)|&\le\big(C_0\norm{m_\phi}p+|\alpha|\norm{\nabla\cdot\eIv{\xi}{}}p\big)\norm{\psi_h}p
+\big(\norm{\eIv{\xi}{}}\Gamma+\norm{\boldsymbol z_u}\Gamma+L^{-1}\norm{z_\phi}\Gamma\big)\norm{\psi_h}\Gamma.
\end{aligned}
\end{equation*}
It remains to prove \eqref{eq:ellbound}.
% REV8-DEL: "Lemma~\ref{lem:projection} implies" (before eq:ellbound; eq:ellbound is now in the statement)
Indeed, by \eqref{eq:F3}--\eqref{eq:F3alt} and Korn's inequality,
$|F_3(t;\bz)|\le C\norm{\eIs{\phi}{}(t)}{H^1(\Omega_p)}\norm{\bz}{a_e}$ for all $\bz\in\Vp$,
and $\partial_tF_3$ satisfies the same bound with $\partial_t\eIs{\phi}{}$ in place of $\eIs{\phi}{}$;
\eqref{eq:approximation} with $j=0,1$ then gives \eqref{eq:ellbound}.
% REV-DEL (item 6): previous justification via the volume part and the interface part of (eq:elldef).
\end{proof}

The bound for $F_2$ thus  requires only a bulk $L^2$ norm and
a tangential trace norm. A direct estimate of   $F_3(\ehv{\xi}{})$  would instead
require spatial control of the structure-velocity error that is not
provided by the energy estimate. We   therefore estimate  its time sum .
% REV8-DEL: "This explains why the bound for $J$ requires" and "We treat its time sum next."
 The key observation is that   $\ehv{\xi}{}=\dt \ehv{\eta}{}+\boldsymbol\kappa$ .
We can therefore transfer the time difference from the displacement error
to the known load   $F_3$ . This replaces a velocity estimate by a
displacement-energy estimate.
 % REV8-DEL: paragraph headings "Size and time variation of the load.", "A discrete integration-by-parts
% formula.", and "Application to the mechanical load."
 The  following lemma transfers the BDF2 difference from a sequence to its
coefficient, leaving only endpoint values and first differences of that
coefficient.

\begin{lemma}[BDF2 summation by parts]\label{lem:sbp}
Let $y^j$ belong to a normed space $Y$ and $\ell^j\in Y'$. For $M\ge2$,
  \begin{equation}
\left|\Delta t\sum_{j=2}^M\ell^j(\dt y^j)\right|
\le C\left(\max_{2\le j\le M}\|\ell^j\|_{Y'}
+\sum_{j=2}^{M-1}\|\ell^{j+1}-\ell^j\|_{Y'}\right)
\max_{0\le j\le M}\|y^j\|_Y,
\label{eq:sbpbound}
\end{equation}  
with an absolute constant independent of $\Delta t$ and $M$.
\end{lemma}
\begin{proof}
 % REV-DEL (item 5): the expansion (eq:sbpexact) valid for M>=3, the interior-coefficient identity,
% and the separate treatment of M=2.
 For $j\ge2$,
$2\Delta t\,\dt y^j=(3y^j-y^{j-1})-(3y^{j-1}-y^{j-2})$.
Summation by parts therefore gives, for every $M\ge2$,
\begin{equation}
\begin{aligned}
\Delta t\sum_{j=2}^{M}\ell^j(\dt y^j)
={}&\frac12\ell^M(3y^M-y^{M-1})
-\frac12\ell^2(3y^1-y^0)\\
&-\frac12\sum_{j=2}^{M-1}
(\ell^{j+1}-\ell^j)(3y^j-y^{j-1}),
\end{aligned}
\label{eq:sbpexact}
\end{equation}  
  where the last sum is empty for  $M=2$ . Since
$\|3y^j-y^{j-1}\|_Y\le4\max_{0\le i\le M}\|y^i\|_Y$, the triangle inequality
yields \eqref{eq:sbpbound} with $C=4$ .
\end{proof}

  \begin{lemma}[Control of the pressure-projection load]\label{lem:loadcontrol}
%\eightmark{from former \S5.5}
For $2\le M\le N$ and every $\epsilon>0$,
\begin{equation}
\left|\Delta t\sum_{j=2}^MF_3(t_j;\ehv{\xi}{j})\right|
\le\epsilon\max_{0\le j\le M}\norm{\ehv{\eta}{j}}{a_e}^2
+C_{T,\epsilon}\big(h^{2k}+\Delta t^4\big)\mathcal M^2.
\label{eq:ellcontrol}
\end{equation}  
 where $C_{T,\epsilon}$ is independent of $h$, $\Delta t$, and $M$.
\end{lemma}
\begin{proof}
Using \eqref{eq:error_kinematic_relation}, write
\begin{equation}
\Delta t\sum_{j=2}^MF_3(t_j;\ehv{\xi}{j})
=\Delta t\sum_{j=2}^MF_3(t_j;\dt \ehv{\eta}{j})
+\Delta t\sum_{j=2}^MF_3(t_j;\boldsymbol\kappa^j).
\label{eq:ellsplit}
\end{equation}
 The fundamental theorem of calculus gives
  \[
\sum_{j=2}^{M-1}\|F_3(t_{j+1})-F_3(t_j)\|_{a_e'}
\le\int_{t_2}^{t_M}\|\partial_tF_3(t)\|_{a_e'}\,dt.
\]  
Thus the variation is bounded by \eqref{eq:ellbound}, independently of
the number of time steps.  %Set $Y_M=\max_{0\le j\le M}\norm{\ehv{\eta}{j}}{a_e}$. 
Applying Lemma~\ref{lem:sbp} to the first sum in \eqref{eq:ellsplit},   with $Y=\Vp$ normed by $\norm{\cdot}{a_e}$, $\ell^j=F_3(t_j)$, and $y^j=\ehv{\eta}{j}$, and  using \eqref{eq:kappabound} for the second, gives
  \[
\left|\Delta t\sum_{j=2}^MF_3(t_j;\ehv{\xi}{j})\right|
\le C_T h^k\mathcal M \max_{0\le j\le M}\norm{\ehv{\eta}{j}}{a_e}
+C_T h^k\Delta t^2\mathcal M^2.
\]  
Young's inequality   and $h^k\Delta t^2\le\frac12(h^{2k}+\Delta t^4)$ give \eqref{eq:ellcontrol}.
% REV8-DEL: "therefore gives, for every fixed $\epsilon>0$," and the display (eq:ellcontrol), now in the statement.
% REV-DEL (item 7, Sec. 5.5): "This treatment uses the time variation of the pressure-projection
% functional and the displacement energy, and does not introduce an inverse trace on normal
% structure velocity." (repeats the opening of Sec. 5; margin flag at the next paragraph)
\end{proof}

\subsection{Energy estimate and convergence theorem}
  \label{subsec:errortheorem}
Let $\E_{h}^j$, $\N_{h}^j$, and $\Diss_{*,h}^j$  be \eqref{eq:discreteenergy}, \eqref{eq:numericaldiss}, and \eqref{eq:Dstar} with all variables replaced by their discrete errors.

 \begin{theorem}[Bulk energy error]\label{thm:error}
%\eightmark{stated before its proof}
Assume the model, mesh, trace, regularity, and finite element approximation hypotheses above. Suppose the initialization satisfies
  \begin{equation}
\E_{h}^1\le C_{\rm init}(h^{2k}+\Delta t^4).
\label{eq:initialization}
\end{equation}  
  If $\Delta t\le C_*h$, then
\begin{align}
\max_{1\le j\le N}\Big(&\norm{\ehv{u}{j}}f^2+\norm{\ehv{\xi}{j}}p^2+
\norm{\ehs{\phi}{j}}p^2+\norm{\ehv{\eta}{j}}{a_e}^2\Big)\notag\\
&+\Delta t\sum_{j=2}^N\left(
\norm{\ehv{u}{j}}{H^1(\Omega_f)}^2+
\norm{\ehs{\phi}{j}}{H^1(\Omega_p)}^2+
\norm{\bP\ehv{\xi}{j}}\Gamma^2\right)
\le C_T(h^{2k}+\Delta t^4).
\label{eq:errorprojected}
\end{align}  
The same bound holds for the total errors $\bu-\bu_h$, $\bxi-\bxi_h$, $\phi-\phi_h$, and $\be-\be_h$, with the indicated norms. In particular,
\begin{equation}
\max_{1\le j\le N}\left(
\norm{\bu(t_j)-\bu_h^j}f+
\norm{\bxi(t_j)-\bxi_h^j}p+
\norm{\phi(t_j)-\phi_h^j}p+
\norm{\be(t_j)-\be_h^j}{H^1(\Omega_p)}\right)
\le C_T(h^k+\Delta t^2).
\label{eq:errortotal}
\end{equation}
The constant can depend on $T$, fixed coefficients, projection constants, $\mathcal M$, and initialization quality, but is independent of $h,\Delta t$.
\end{theorem}
\begin{proof}
  We now combine the preceding bounds with the interface stability
estimates. %\revdel{} (REV11: space restored before the comment)
The argument is organized into the energy identity, residual
control, and a summation step in which the maximum of the error energy is absorbed; no discrete Gronwall inequality is used .

\emph{ Step 1: The error energy identity .}
 % REV8-DEL: "Let E_h^j, N_h^j, and D_{*,h}^j be ... replaced by their discrete errors." (moved before Theorem thm:error)
 Testing \eqref{eq:errorfluid}--\eqref{eq:errorbiot} with $(\ehv{u}{n+1},\ehs{p}{n+1},\ehv{\xi}{n+1},\ehs{\phi}{n+1})$ gives
\begin{align}
\frac{\E_{h}^{n+1}-\E_{h}^n}{\Delta t}+\N_{h}^{n+1}
+2\mu_f\norm{\bD(\ehv{u}{n+1})}f^2
+\norm{K^{1/2}\nabla\ehs{\phi}{n+1}}p^2
&+\gamma\norm{\bP(\ehv{u}{n+1}-\ehv{\xi}{n+1})}\Gamma^2\notag\\
=\R_\Gamma^{n+1}(e^h)
+F_1^{n+1}(\ehv{u}{n+1})+F_2^{n+1}(\ehv{\xi}{n+1})
+F_4^{n+1}(\ehs{\phi}{n+1})
&+F_3(t_{n+1};\ehv{\xi}{n+1})
-a_e(\ehv{\eta}{n+1},\boldsymbol\kappa^{n+1}).
\label{eq:errorenergy}
\end{align}  
  The last term follows from $a_e(\ehv{\eta}{},\ehv{\xi}{})=a_e(\ehv{\eta}{},\dt \ehv{\eta}{})+a_e(\ehv{\eta}{},\boldsymbol\kappa)$; there is no additional kinematic testing term .

% REV4-DEL: former Step 2 (bound of the residuals by (1/2)D + C E^{n+1} + ..., eq. (errorresidualbound))
% and former Step 3 (summed inequality with C Delta t sum E^j, eq. (errorpresum)); replaced below.
\emph{ Step 2: The per-step error inequality .}
  All arguments of $\R_\Gamma^{n+1}(e^h)$ are finite element functions, so
Steps~2 and~3 of the proof of Theorem~\ref{thm:stability} apply verbatim under
$\Delta t\le C_*h$: $\R_\Gamma^{n+1}(e^h)$ satisfies \eqref{eq:Rbound}, and
\eqref{eq:tangentialcontrol} holds for $(\ehv{u}{n+1},\ehv{\xi}{n+1})$. Hence
\eqref{eq:errorenergy} gives
\[
\begin{aligned}
\frac{\E_h^{n+1}-\E_h^n}{\Delta t}+\frac12\N_h^{n+1}+\Diss_{*,h}^{n+1}
\le{}&F_1^{n+1}(\ehv{u}{n+1})+F_2^{n+1}(\ehv{\xi}{n+1})+F_4^{n+1}(\ehs{\phi}{n+1})\\
&+F_3(t_{n+1};\ehv{\xi}{n+1})-a_e(\ehv{\eta}{n+1},\boldsymbol\kappa^{n+1}).
\end{aligned}
\]
By \eqref{eq:FGJbounds}, Korn's and Poincar\'e's inequalities, and
Young's inequality, using $(h^k+\Delta t^2)^2\le2(h^{2k}+\Delta t^4)$,
\[
\begin{aligned}% REV8: split into two lines
|F_1^{n+1}(\ehv{u}{n+1})|&\le\frac{\mu_f}{2}\norm{\bD(\ehv{u}{n+1})}f^2+C(h^{2k}+\Delta t^4)\mathcal M^2,\\
|F_4^{n+1}(\ehs{\phi}{n+1})|&\le\frac14\norm{K^{1/2}\nabla\ehs{\phi}{n+1}}p^2+C(h^{2k}+\Delta t^4)\mathcal M^2.
\end{aligned}
\]
We split $F_2^{n+1}(\ehv{\xi}{n+1})$ into its bulk part
$\rho_p(\boldsymbol m_\xi^{n+1},\ehv{\xi}{n+1})_p$ and its tangential part; by
\eqref{eq:projectiontrace} and \eqref{eq:defectbounds},
\[
\gamma\left|\dual{\bP(\eIv{\xi}{n+1}-\boldsymbol z_u^{n+1})}{\bP\ehv{\xi}{n+1}}\right|
\le\frac{q_\Gamma}{4}\norm{\bP\ehv{\xi}{n+1}}\Gamma^2+C(h^{2k}+\Delta t^4)\mathcal M^2.
\]
These three absorptions use one half of $\Diss_{*,h}^{n+1}$. Hence, for $n\ge1$,
% REV8-DEL: "Since $\chi^2\le2(h^{2k}+\Delta t^4)$, we obtain"
\begin{equation}
\begin{aligned}
\frac{\E_h^{n+1}-\E_h^n}{\Delta t}+\frac12\N_h^{n+1}+\frac12\Diss_{*,h}^{n+1}
\le{}&\rho_p(\boldsymbol m_\xi^{n+1},\ehv{\xi}{n+1})_p
-a_e(\ehv{\eta}{n+1},\boldsymbol\kappa^{n+1})\\
&+F_3(t_{n+1};\ehv{\xi}{n+1})
+C_1\big(h^{2k}+\Delta t^4\big)\mathcal M^2.
\end{aligned}
\label{eq:errorresidualbound}
\end{equation}
The first two terms on the right are tested with the bulk structure-velocity
and displacement errors, which are not directly controlled by the dissipation. We do not bound them by the
current energy, which would lead to a Gronwall argument; instead, we sum them in
time and use the maximum of the error energy.

\emph{Step 3: Summation in time.}
Let $2\le m\le M\le N$. Multiplying \eqref{eq:errorresidualbound} by
$\Delta t$ and summing over $n=1,\ldots,m-1$ gives, since $(m-1)\Delta t\le T$,
\begin{equation}
\begin{aligned}
\E_h^m+\frac{\Delta t}{2}\sum_{j=2}^m\big(\N_h^j+\Diss_{*,h}^j\big)
\le{}&\E_h^1+C_1T\big(h^{2k}+\Delta t^4\big)\mathcal M^2
+\left|\Delta t\sum_{j=2}^mF_3(t_j;\ehv{\xi}{j})\right|\\
&+\Delta t\sum_{j=2}^m\Big(\rho_p\norm{\boldsymbol m_\xi^j}p\norm{\ehv{\xi}{j}}p
+\norm{\ehv{\eta}{j}}{a_e}\norm{\boldsymbol\kappa^j}{a_e}\Big).
\end{aligned}
\label{eq:errorpresum}
\end{equation}
Let $S_M=\max_{1\le j\le M}\E_h^j$. By the definition of the BDF2 energy,
$\rho_p\norm{\ehv{\xi}{j}}p^2\le4\E_h^j$ and $\norm{\ehv{\eta}{j}}{a_e}^2\le4\E_h^j$
for $j\ge1$. Moreover, the identity
$\ehv{\eta}{0}=2\ehv{\eta}{1}-(2\ehv{\eta}{1}-\ehv{\eta}{0})$ gives
$\norm{\ehv{\eta}{0}}{a_e}^2\le5\big(\norm{\ehv{\eta}{1}}{a_e}^2+\norm{2\ehv{\eta}{1}-\ehv{\eta}{0}}{a_e}^2\big)\le20\E_h^1$.
Hence
\begin{equation}
\max_{2\le j\le M}\norm{\ehv{\xi}{j}}p\le2\rho_p^{-1/2}S_M^{1/2},
\qquad
\max_{0\le j\le M}\norm{\ehv{\eta}{j}}{a_e}^2\le20S_M.
\label{eq:maxcontrol}
\end{equation}
We bound the last three terms of \eqref{eq:errorpresum}. By
\eqref{eq:consistencybounds}, \eqref{eq:maxcontrol}, and Young's inequality,
\[
\Delta t\sum_{j=2}^m\rho_p\norm{\boldsymbol m_\xi^j}p\norm{\ehv{\xi}{j}}p
\le CT(h^k+\Delta t^2)\mathcal M\,S_M^{1/2}
\le\frac18S_M+CT^2\big(h^{2k}+\Delta t^4\big)\mathcal M^2.
\]
By the continuity of $a_e$, \eqref{eq:kappabound}, and \eqref{eq:maxcontrol},
\[
\Delta t\sum_{j=2}^m\norm{\ehv{\eta}{j}}{a_e}\norm{\boldsymbol\kappa^j}{a_e}
\le CT\Delta t^2\mathcal M\,S_M^{1/2}
\le\frac18S_M+CT^2\Delta t^4\mathcal M^2.
\]
Finally, \eqref{eq:ellcontrol} with $M$ replaced by $m$  and   $\epsilon=1/80$, together
with \eqref{eq:maxcontrol}, bounds the $F_3$-term by
$\frac14S_M+C_{T}(h^{2k}+\Delta t^4)\mathcal M^2$. Combining these bounds with
\eqref{eq:errorpresum} yields, for every $2\le m\le M$,
\begin{equation}
\E_h^m+\frac{\Delta t}{2}\sum_{j=2}^m\big(\N_h^j+\Diss_{*,h}^j\big)
\le\E_h^1+\frac12S_M+C_T\big(h^{2k}+\Delta t^4\big)\mathcal M^2,
\label{eq:errormaxbound}
\end{equation}
where $C_T$ is independent of $h$, $\Delta t$, $m$, and $M$.

% REV4-DEL: former proof steps "Control of the starting levels" (moved to Step 3),
% "Absorption of the displacement maximum", and "Discrete Gronwall".
\emph{Step 4: Absorption of the maximum.}
Write $K_T=C_T(h^{2k}+\Delta t^4)\mathcal M^2$ for the last term in
\eqref{eq:errormaxbound}. Inequality \eqref{eq:errormaxbound} also holds for $m=1$,
where the sum is empty. Dropping the nonnegative dissipation and taking the maximum
over $1\le m\le M$ gives $S_M\le\E_h^1+\frac12S_M+K_T$, that is,
$S_M\le2\E_h^1+2K_T$. Inserting this into \eqref{eq:errormaxbound} with $m=M$ gives
\[
\E_h^M+\frac{\Delta t}{2}\sum_{j=2}^M\big(\N_h^j+\Diss_{*,h}^j\big)
\le2\E_h^1+2K_T,\qquad 2\le M\le N.
\]

\emph{Step 5: Dissipation and total errors.}
The bound for $S_N$ controls the maximum in \eqref{eq:errorprojected}, and the last display controls the summed numerical and physical dissipation. Korn and Poincar\'{e} inequalities convert the fluid and  scalar parts of   $\Diss_{*,h}$ into the stated  $H^1$   norms. Finally, use \eqref{eq:initialization}, the approximation assumptions,  the   triangle inequality, and  $h\le1$ to obtain the total-error estimate.
\end{proof}

\section{Comparison with Robin-Robin splitting scheme}
\label{sec:comparison}

 % REV-DEL (item 7, Sec. 6): the opening paragraph restating the difference between the two
% schemes (already given in Sec. 3.3, eqs. (eq:RRnormal)-(eq:normaldifference)).
In this section we identify the residual that leads to the sufficient restriction $\Delta t\lesssim h^2$ in the analysis of the Robin--Robin splitting in \cite{WangLeeCanic2026SecondOrder}.
% REV10-DEL: "we explain why the Robin-Robin splitting in [1] has a difficulty to avoid the CFL-type condition Delta t <= Ch^2."
The mechanical normal condition in \cite{WangLeeCanic2026SecondOrder} has the form
\begin{equation}
\bn_p\cdot\bs_p^{n+1}\bn_p+\phi^{n+1}
+\beta\bxi^{n+1}\cdot\bn_p
=\beta\bxi^{*,n+1}\cdot\bn_p,
\label{eq:RRnormal}
\end{equation}
 with a fixed Robin parameter $\beta>0$. Compared with \eqref{eq:schememech}, its weak equation contains two additional terms, namely
\begin{equation}
\beta\dual{\bxi_h^{n+1}\cdot\bn_p}{\bz_h\cdot\bn_p}
\quad\text{on the left},\qquad
\beta\dual{\bxi_h^{*,n+1}\cdot\bn_p}{\bz_h\cdot\bn_p}
\quad\text{on the right}.
\label{eq:deletedpair}
\end{equation}
In residual form, the difference is
\begin{equation}
\beta\dual{(\bxi_h^{n+1}-2\bxi_h^n+\bxi_h^{n-1})\cdot\bn_p}{\bz_h\cdot\bn_p}.
\label{eq:normaldifference}
\end{equation}
For a smooth exact solution this is $O(\Delta t^2)$, so removing it preserves formal second-order consistency. It nevertheless changes the discrete method.
Thus the two schemes differ only by the normal structure-velocity residual \eqref{eq:normaldifference}, with $\br_\xi^{n+1}=\bxi_h^{n+1}-2\bxi_h^n+\bxi_h^{n-1}$.
The additional normal-velocity term contributes
\begin{equation}
-\beta\dual{\br_\xi^{n+1}\cdot\bn_p}
{\bxi_h^{n+1}\cdot\bn_p}
\label{eq:RRbadterm}
\end{equation}
to the energy identity. If both interface traces are estimated through
bulk $L^2$ norms, each inverse trace inequality contributes a factor
$h^{-1/2}$. Cauchy--Schwarz and Young's inequality then give
\begin{align}
\beta\left|\dual{\br_\xi\cdot\bn_p}{\bxi_h\cdot\bn_p}\right|
&\le C\beta h^{-1}\norm{\br_\xi}p\norm{\bxi_h}p\notag \\
&\le \frac{\rho_p}{32\Delta t}\norm{\br_\xi}p^2
+C\frac{\beta^2\Delta t}{\rho_ph^2}\norm{\bxi_h}p^2.
\label{eq:RRinverse}
\end{align}
The first term is absorbed by the BDF2 numerical dissipation.
The second is proportional to the current structure kinetic energy,
with a coefficient of order $\Delta t/h^2$.

To track this contribution, let $\E_{\mathrm{RR}}^j$,
$\N_{\mathrm{RR}}^j$, and $\Diss_{*,\mathrm{RR}}^j$ be the quantities
in \eqref{eq:discreteenergy}, \eqref{eq:numericaldiss}, and
\eqref{eq:Dstar}, evaluated on the Robin--Robin variant retaining
\eqref{eq:deletedpair}.
Including \eqref{eq:RRinverse}, the total coefficient assigned to
$\norm{\br_\xi^{n+1}}p^2$ is
$2\rho_p/(32\Delta t)$, which still fits within
$\N_{\mathrm{RR}}^{n+1}/2$.
Since
\[
\E_{\mathrm{RR}}^{n+1}
\ge\frac{\rho_p}{4}\norm{\bxi_h^{n+1}}p^2,
\]
the remaining term yields, for zero sources and $\Delta t\le C_*h$,
\begin{equation}
\frac{\E_{\mathrm{RR}}^{n+1}-\E_{\mathrm{RR}}^n}{\Delta t}
+\frac12\N_{\mathrm{RR}}^{n+1}+\Diss_{*,\mathrm{RR}}^{n+1}
\le a_{\mathrm{RR}}\frac{\Delta t}{h^2}\E_{\mathrm{RR}}^{n+1}.
\label{eq:RRgrowth}
\end{equation}
Here $a_{\mathrm{RR}}>0$ is independent of $h$ and $\Delta t$ and
contains the fixed factor $\beta^2/\rho_p^2$.

\paragraph{Accumulation over time.}
Set $z=a_{\mathrm{RR}}\Delta t^2/h^2$ and suppose $z\le1/2$.
Multiplying \eqref{eq:RRgrowth} by $\Delta t$ gives
\begin{equation}
(1-z)\E_{\mathrm{RR}}^{n+1}
+\Delta t\left(\frac12\N_{\mathrm{RR}}^{n+1}
+\Diss_{*,\mathrm{RR}}^{n+1}\right)
\le\E_{\mathrm{RR}}^n.
\label{eq:RRrecursion}
\end{equation}
Thus each step can increase the energy by a factor $(1-z)^{-1}$.
Iteration, retaining the accumulated nonnegative dissipation, gives
\begin{align}
&\E_{\mathrm{RR}}^M
+\Delta t\sum_{j=2}^{M}
\left(\frac12\N_{\mathrm{RR}}^j+\Diss_{*,\mathrm{RR}}^j\right)
\notag\\
&\qquad\le(1-z)^{-(M-1)}\E_{\mathrm{RR}}^1
\le\exp\!\left(2a_{\mathrm{RR}}T\frac{\Delta t}{h^2}\right)
\E_{\mathrm{RR}}^1,
\qquad t_M\le T.
\label{eq:RRgronwall}
\end{align}
The last inequality uses $(1-z)^{-1}\le e^{2z}$ and
$(M-1)\Delta t\le T$.
The exponential is bounded uniformly under refinement if
$\Delta t/h^2$ remains bounded. This is why this argument gives a
parabolic time-step restriction for the earlier scheme.

More precisely, because $h\le1$, the choice
\[
\Delta t\le c_{\mathrm{RR}}h^2,
\qquad
0<c_{\mathrm{RR}}\le
\min\left\{C_*,(2a_{\mathrm{RR}})^{-1/2}\right\}
\]
ensures both $\Delta t\le C_*h$ and $z\le1/2$.
It also bounds the exponential in \eqref{eq:RRgronwall} by
$\exp(2a_{\mathrm{RR}}c_{\mathrm{RR}}T)$, independently of the mesh
and time step.

In the Neumann-Robin splitting scheme, \eqref{eq:RRbadterm} does not appear,
so it removes the growth contribution on the
right of \eqref{eq:RRgrowth}. The remaining interface residuals are
controlled under $\Delta t\le C_*h$, giving the closed energy estimate
of Theorem~\ref{thm:stability}.
 %The improvement is therefore in the sufficient time-step condition
%proved by this energy argument: $\Delta t\lesssim h$ instead of
%$\Delta t\lesssim h^2$.
%It does not establish that the parabolic restriction is necessary for
%the earlier scheme, or that either bound is a sharp stability threshold.

\begin{table}[t]
\centering
\caption{Comparison with the second-order Robin--Robin scheme in
\cite{WangLeeCanic2026SecondOrder}. The time-step conditions are
sufficient analytical bounds for fixed coefficients.}
\label{tab:comparison}
\small
\begin{tabularx}{\textwidth}{@{}>{\raggedright\arraybackslash}p{0.29\textwidth}>{\raggedright\arraybackslash}X>{\raggedright\arraybackslash}X@{}}
\toprule
Feature &  Robin--Robin scheme in \cite{WangLeeCanic2026SecondOrder} & Present Neumann--Robin scheme\\
\midrule
Mechanical normal condition
& Traction with an auxiliary normal-velocity Robin term
& Physical Neumann-type traction with current pore pressure\\
Normal structure-velocity predictor in Biot traction
& Present & Absent\\
Fluid normal and Darcy conditions & Robin & Robin\\
Tangential coupling & Robin/BJS & Robin/BJS\\
Time discretization & BDF2 with EX2 & BDF2 with EX2\\
Interface data from the other block & Previous levels & Previous levels\\
Subdomain solves & Independent and parallel & Independent and parallel\\
Sufficient time-step condition
& $\Delta t\lesssim h^2$ & $\Delta t\lesssim h$\\
\bottomrule
\end{tabularx}
\end{table}

\section{Numerical experiments}
\label{sec:verification}

\subsection{Test 1: spatial and temporal convergence study}
\label{test1}

We assess the accuracy of the second-order Neumann--Robin splitting
scheme using a polynomial manufactured solution. The exact solutions determine the volume sources and satisfy all physical interface conditions without additional interface loads. Separate spatial and temporal refinement studies are performed to examine the corresponding convergence.

Consider the following domains for fluid, poroelastic structure and the interface
\[
\Omega_f=(0,1)\times(0,1),\qquad
\Omega_p=(0,1)\times(-1,0),\qquad
\Gamma=(0,1)\times\{0\},
\]
with outward interface normals
$\bn_f=(0,-1)^T$ and $\bn_p=(0,1)^T$.
We set the final time $T=1$ and choose the parameters
$\rho_f=\rho_p=\mu_f=\mu_p=\lambda_p
=\alpha=C_0=\gamma=L=1, K=\boldsymbol I$.

Define the spatial polynomials
\begin{align}
F(x)&=x^2(1-x)^2,
\quad A(y)=(1+y)^2(1-2y),
\quad B(y)=y^2(1+y),\notag\\
C(y)&=y-3y^3+2y^4,
\quad D(y)=y(1+y)^2.\notag
\label{eq:MMSdefinitions}
\end{align}
Using the stream function
\begin{equation}
\Psi(t,x,y)=2tF(x)C(y),
\label{eq:MMSstream}
\end{equation}
we prescribe the exact fields
\begin{equation}
\left\{
\begin{aligned}
\bu(t,x,y)
&=2t\big(F(x)C'(y),-F'(x)C(y)\big)^T,\\
p_f(t,x,y)
&=t^2F(x)-4tF'(x),\\
\be(t,x,y)
&=t^2\left(F(x)A(y),-\frac13F'(x)D(y)\right)^T,\\
\bxi(t,x,y)
&=2t\left(F(x)A(y),-\frac13F'(x)D(y)\right)^T,\\
\phi(t,x,y)
&=F(x)\big(B(y)+t^2A(y)\big).
\end{aligned}
\right.
\label{eq:MMSfields}
\end{equation}
In particular, the stream-function construction gives
$\bu=(\partial_y\Psi,-\partial_x\Psi)^T$ and
$\nabla\cdot\bu=0$. The kinematic relation $\partial_t\be=\bxi$
also holds identically.
Moreover, the manufactured solutions also ensure the homogeneous exterior conditions as follows
\[
\bu=0\quad\text{on }\Sigma_f,\qquad
\be=\bxi=0,\quad \phi=0
\quad\text{on }\Sigma_p.
\]
No Dirichlet condition is imposed on the fluid pressure, and its level
is determined by the interface traction condition. At the interface $\Gamma$, the interface traces satisfy that
\begin{align}
\bu|_\Gamma=\bxi|_\Gamma=(2tF(x),0)^T,\quad
\be|_\Gamma=(t^2F(x),0)^T,\quad
\phi|_\Gamma=t^2F(x),\quad
\partial_y\phi|_\Gamma&=0.
\label{eq:MMStraces}
\end{align}
We note that the manufactured solution has polynomial temporal dependence of degree at most two. Consequently, the BDF2 approximation is exact when applied to the corresponding exact solution values. The second-order extrapolation of the pore pressure, however, introduces a nonzero temporal consistency defect of order $(\Delta t^2)$. Moreover, the exact normal structure velocity vanishes on the interface, namely $\boldsymbol{\xi}\cdot\boldsymbol{n}_p=0$. Thus, the present example provides a convergence test for the coupled numerical method, but does not independently test the normal structure-velocity residual removed by the proposed Neumann–Robin modification. 

The volume sources are obtained by substituting the exact fields
into the continuous equations:
\begin{equation}
\left\{
\begin{aligned}
\boldsymbol f_f
&=\rho_f\partial_t\bu
  -\nabla\cdot\bs_f(\bu,p_f),\notag\\
\boldsymbol f_p
&=\rho_p\partial_t\bxi
  -\nabla\cdot\bs_p(\be,\phi),\notag\\
g_p
&=C_0\partial_t\phi+\alpha\nabla\cdot\bxi
  -\nabla\cdot(K\nabla\phi).
\end{aligned}
\right.
\label{eq:MMSsources}
\end{equation}
And the initial data are
\[
\bu(0)=0,\qquad p_f(0)=0,\qquad
\be(0)=\bxi(0)=0,\qquad
\phi(0,x,y)=F(x)B(y).
\]
In particular, the initial pore pressure is nonzero in the interior
of $\Omega_p$. The two BDF2 starting levels, at $t_0=0$ and
$t_1=\Delta t$, are prescribed by projecting the exact fields.
We use a joint Robin--Stokes projection for the fluid pair,
a common elasticity Ritz projection for displacement and structure
velocity, and the Darcy Ritz projection defined in \eqref{eq:ritzmaps} for pore
pressure.

In both refinement studies, we use continuous piecewise quadratic
($P_2$) elements for the fluid velocity, poroelastic displacement,
structure velocity, and pore pressure, and continuous piecewise
linear ($P_1$) elements for the fluid pressure.
Both subdomains are discretized using uniform triangular meshes
with $n$ intervals in each coordinate direction. We report the
mesh spacing as $\Delta x=1/n$.

Let $t_j=j\Delta t$ and $N=T/\Delta t$. The displacement error is measured
in the elastic energy norm:
\begin{equation}
E_\eta
=\max_{2\le j\le N}
\left(
2\mu_p
\|\bD(\be(t_j)-\be_h^j)\|_{L^2(\Omega_p)}^2
+\lambda_p
\|\nabla\cdot(\be(t_j)-\be_h^j)\|_{L^2(\Omega_p)}^2
\right)^{1/2}.
\label{eq:numericaletaerror}
\end{equation}
Under the homogeneous displacement boundary conditions, this norm
is equivalent to the $H^1(\Omega_p)$ norm.
For the remaining variables, we use maximum-in-time $L^2$ errors:
\[
E_z=\max_{2\le j\le N}
\|z(t_j)-z_h^j\|_{L^2(\Omega_*)},
\qquad z\in\{\bu,p_f,\bxi,\phi\},
\]
where $\Omega_*=\Omega_f$ for $\bu,p_f$ and $\Omega_*=\Omega_p$
for $\bxi,\phi$. All maxima are taken over the computed BDF2
time levels.

For spatial refinement, we chooce $\Delta t=10^{-3}$ and $T=1$ and use
$n=4,8,16,32,64,128$. The observed rates are calculated as
\[
\mathrm{rate}
=\frac{\log(E_{\Delta x}/E_{\Delta x/2})}{\log 2}.
\]

\begin{table}[htbp]
\centering
\caption{Spatial convergence of the Neumann--Robin scheme for
$T=1$ and fixed $\Delta t=10^{-3}$.
The displacement error $E_\eta$ is measured in the
$H^1$-equivalent elastic energy norm, and all other errors are
measured in $L^2$. All errors are maxima over the computed
BDF2 time levels.}
\label{tab:spatial-convergence}

\medskip
\textbf{(a) Fluid variables}\par\smallskip

\begingroup
\scriptsize
\renewcommand{\arraystretch}{1.1}
\setlength{\tabcolsep}{2pt}
\begin{tabular*}{\linewidth}{@{\extracolsep{\fill}}ccrrrr@{}}
\toprule
$n$ & $\Delta x$ & $E_u$ & Rate & $E_{p_f}$ & Rate\\
\midrule
4 & $1/4$
& $2.79\times10^{-3}$ & --
& $7.22\times10^{-2}$ & --\\

8 & $1/8$
& $3.14\times10^{-4}$ & 3.15
& $1.66\times10^{-2}$ & 2.12\\

16 & $1/16$
& $3.63\times10^{-5}$ & 3.11
& $4.06\times10^{-3}$ & 2.03\\

32 & $1/32$
& $4.40\times10^{-6}$ & 3.05
& $1.01\times10^{-3}$ & 2.01\\

64 & $1/64$
& $5.43\times10^{-7}$ & 3.02
& $2.52\times10^{-4}$ & 2.00\\

128 & $1/128$
& $6.76\times10^{-8}$ & 3.01
& $6.31\times10^{-5}$ & 2.00\\
\bottomrule
\end{tabular*}
\endgroup

\medskip
\textbf{(b) Poroelastic variables}\par\smallskip

\begingroup
\scriptsize
\renewcommand{\arraystretch}{1.1}
\setlength{\tabcolsep}{2pt}
\begin{tabular*}{\linewidth}{@{\extracolsep{\fill}}ccrrrrrr@{}}
\toprule
$n$ & $\Delta x$
& $E_\eta$ & Rate
& $E_\xi$ & Rate
& $E_\phi$ & Rate\\
\midrule
4 & $1/4$
& $1.84\times10^{-2}$ & --
& $1.09\times10^{-3}$ & --
& $4.21\times10^{-4}$ & --\\

8 & $1/8$
& $4.98\times10^{-3}$ & 1.88
& $1.28\times10^{-4}$ & 3.09
& $5.50\times10^{-5}$ & 2.93\\

16 & $1/16$
& $1.28\times10^{-3}$ & 1.96
& $1.51\times10^{-5}$ & 3.09
& $6.98\times10^{-6}$ & 2.98\\

32 & $1/32$
& $3.25\times10^{-4}$ & 1.98
& $1.84\times10^{-6}$ & 3.04
& $8.78\times10^{-7}$ & 2.99\\

64 & $1/64$
& $8.17\times10^{-5}$ & 1.99
& $2.27\times10^{-7}$ & 3.01
& $1.10\times10^{-7}$ & 2.99\\

128 & $1/128$
& $2.05\times10^{-5}$ & 2.00
& $2.85\times10^{-8}$ & 3.00
& $1.69\times10^{-8}$ & 2.71\\
\bottomrule
\end{tabular*}
\endgroup
\end{table}
Table~\ref{tab:spatial-convergence} shows approximately second-order spatial convergence for the displacement in the elastic energy norm and the fluid pressure in the $L^2$ norm. The fluid velocity, structure velocity, and pore pressure exhibit approximately third-order convergence in the $L^2$ norm over most of the refinement range. These observed rates are higher than the $O(h^k)$ spatial convergence rate guaranteed by the bulk energy-norm error estimate. For the $P_2$ approximation ($k=2$), the theorem establishes an $O(h^2+\Delta t^2)$ bound, whereas the numerical results exhibit approximately $O(h^3)$ spatial convergence for several variables in the $L^2$ norm. Such higher-order behavior is consistent with the approximation properties of quadratic finite elements under additional regularity assumptions. However, establishing an optimal $L^2$ error estimate for the fully discrete splitting scheme would require further analysis of the interface and splitting errors and is beyond the scope of the present work.

For temporal refinement, we fix the mesh at $\Delta x=1/250$ intervals in each coordinate direction, and $T=1$. The time step is set to be  halved successively:
\[
\Delta t\in
\left\{\frac14,\frac18,\frac1{16},\frac1{32},\frac1{64}\right\}.
\]
The observed temporal rates are calculated as
\[
\mathrm{rate}
=\frac{\log(E_\Delta t/E_{\Delta t/2})}{\log 2}.
\]

\begin{table}[htbp]
\centering
\caption{Temporal convergence of the Neumann--Robin (NR, $\beta=0$)
scheme for the old polynomial manufactured solution with
$T=1$ on the fixed mesh $n=250$, $\Delta x=1/250$.
The displacement error $E_\eta$ is measured in the
$H^1$-equivalent elastic energy norm, and all other errors are
measured in $L^2$. All errors are maxima over the computed
BDF2 time levels.}
\label{tab:temporal-convergence}

\medskip
\textbf{(a) Fluid variables}\par\smallskip

\begingroup
\scriptsize
\renewcommand{\arraystretch}{1.1}
\setlength{\tabcolsep}{2pt}
\begin{tabular*}{\linewidth}{@{\extracolsep{\fill}}crrrr@{}}
\toprule
$\Delta t$ & $E_u$ & Rate & $E_{p_f}$ & Rate\\
\midrule
$1/4$
& $7.42\times10^{-5}$ & --
& $7.26\times10^{-3}$ & --\\

$1/8$
& $1.83\times10^{-5}$ & 2.02
& $1.77\times10^{-3}$ & 2.03\\

$1/16$
& $4.48\times10^{-6}$ & 2.03
& $4.29\times10^{-4}$ & 2.05\\

$1/32$
& $1.08\times10^{-6}$ & 2.06
& $1.07\times10^{-4}$ & 2.00\\

$1/64$
& $2.70\times10^{-7}$ & 2.00
& $3.13\times10^{-5}$ & 1.78\\
\bottomrule
\end{tabular*}
\endgroup

\medskip
\textbf{(b) Poroelastic variables}\par\smallskip

\begingroup
\scriptsize
\renewcommand{\arraystretch}{1.1}
\setlength{\tabcolsep}{2pt}
\begin{tabular*}{\linewidth}{@{\extracolsep{\fill}}crrrrrr@{}}
\toprule
$\Delta t$
& $E_\eta$ & Rate
& $E_\xi$ & Rate
& $E_\phi$ & Rate\\
\midrule
$1/4$
& $5.75\times10^{-4}$ & --
& $2.71\times10^{-4}$ & --
& $6.31\times10^{-4}$ & --\\

$1/8$
& $1.91\times10^{-4}$ & 1.59
& $8.95\times10^{-5}$ & 1.60
& $1.49\times10^{-4}$ & 2.08\\

$1/16$
& $5.01\times10^{-5}$ & 1.93
& $2.40\times10^{-5}$ & 1.90
& $3.78\times10^{-5}$ & 1.98\\

$1/32$
& $1.35\times10^{-5}$ & 1.89
& $6.17\times10^{-6}$ & 1.96
& $9.55\times10^{-6}$ & 1.98\\

$1/64$
& $6.17\times10^{-6}$ & 1.13
& $1.56\times10^{-6}$ & 1.99
& $2.41\times10^{-6}$ & 1.99\\
\bottomrule
\end{tabular*}
\endgroup
\end{table}
Table~\ref{tab:temporal-convergence} shows rates close to two
for the fluid velocity, structure velocity, and pore pressure
throughout the finer temporal refinements. On the finest
refinement, the observed convergence rates for
$E_u$, $E_\xi$, and $E_\phi$ are $2.00$, $1.99$, and $1.99$,
respectively, supporting the second-order temporal behavior
of the BDF2--EX2 discretization for these variables.
In contrast, the fluid-pressure convergence rate decreases
to $1.78$, while the displacement convergence rate decreases
to $1.13$, suggesting that spatial discretization errors
may influence the finest temporal refinements.

To examine these reductions variable by variable, we compare
the finest temporal errors with spatial errors estimated from
Table~\ref{tab:spatial-convergence}. Scaling the errors for
$n=128$ to $n=250$ with the approximate observed spatial
orders, two for $\be$ and $p_f$ and three for $\bu$ and
$\phi$, gives
\begin{align*}
E_{\eta,\mathrm{space}}(\Delta x=1/250)
&\approx 2.05\times10^{-5}
\left(\frac{128}{250}\right)^2
\approx 5.37\times10^{-6},\\
E_{p_f,\mathrm{space}}(\Delta x=1/250)
&\approx 6.31\times10^{-5}
\left(\frac{128}{250}\right)^2
\approx 1.65\times10^{-5},\\
E_{u,\mathrm{space}}(\Delta x=1/250)
&\approx 6.76\times10^{-8}
\left(\frac{128}{250}\right)^3
\approx 9.07\times10^{-9},\\
E_{\phi,\mathrm{space}}(\Delta x=1/250)
&\approx 1.69\times10^{-8}
\left(\frac{128}{250}\right)^3
\approx 2.27\times10^{-9}.
\end{align*}
Since Table~\ref{tab:spatial-convergence} was computed with
$\Delta t=10^{-3}$, these values approximate, but do not
isolate, the spatial error. The estimated displacement spatial error, $5.37\times10^{-6}$, is close to the finest temporal error,
$6.17\times10^{-6}$, supporting the interpretation that
spatial error limits the final reduction of the displacement
error. Similarly, the finest fluid-pressure error,
$3.13\times10^{-5}$, is approximately $1.89$ times
the estimated spatial error, suggesting that spatial
discretization contributes to its reduced convergence rate.
In contrast, the finest fluid-velocity error,
$2.70\times10^{-7}$, is approximately $30$ times
the estimated spatial error, while the finest pore-pressure
error, $2.41\times10^{-6}$, is approximately $1060$
times the corresponding spatial-error estimate.
These differences are consistent with the nearly
second-order temporal convergence observed for fluid
velocity and pore pressure, whose errors do not appear
to be limited by the estimated spatial-error levels
over the tested range.

In the temporal study, $\Delta t/\Delta x$ ranges from
$3.91$ to $62.5$. Since the constant $C_*$ in
\eqref{eq:cstar} has not been evaluated for this
configuration, we do not claim that these runs satisfy
the sufficient condition $\Delta t\le C_*h$ of
Theorem~\ref{thm:stability}. Thus, the observed convergence
is empirical and does not test the sharpness of that
condition.
\subsection{Test 2}

We compare the second-order Neumann-Robin (NR) splitting scheme with the
Robin-Robin (RR) variant in \cite{WangLeeCanic2026SecondOrder} using the same nontrivial polynomial manufactured solution. The exact fields satisfy the physical
interface conditions without additional interface loads. Unlike 
Test 1 in \ref{test1}, the present solution has nonzero normal skeleton motion, Darcy filtration, and tangential
slip, as well as nonzero temporal consistency defects. Both schemes are tested with 
the same meshes, forcing terms, starting values, and finite element spaces, and
they differ only in their poroelastic mechanical normal condition.

We consider the same domains and parameter values
as depicted in \ref{test1}.
For the Robin-Robin scheme we introduce a parameter $\beta=1$ to denote the additional normal mechanical
component. Then, the Neumann-Robin scheme corresponds to the same mechanical interface equation with $\beta=0$. Namely, the
normal condition used for the comparative tests is
\[
\bn_p\cdot\bs_p(\be_h^{n+1},\phi_h^{n+1})\bn_p
+\phi_h^{n+1}+\beta\bxi_h^{n+1}\cdot\bn_p
=\beta (2\bxi_h^n-\bxi_h^{n-1})\cdot\bn_p.
\]
While all other variational terms are identical in both runs.

Define the time-dependent polynomials
\begin{equation}
 a(t)=t+t^4,\qquad s(t)=t+t^3,\qquad b(t)=1+t^3,
 \label{eq:newMMS_time}
\end{equation}
and the spatial polynomials
\begin{align}
 F(x)&=256x^4(1-x)^4, \qquad
 A(y)=y(1+y), \qquad
 B(y)=(1+y)^2(1-2y),\notag\\
 C(y)&=1-4y^3+3y^4,\qquad
 D(y)=y-3y^3+2y^4,\qquad
 H(y)=\tfrac12 y^2(1-y)^2.
 \label{eq:newMMS_definitions}
\end{align}
Using the stream function
\begin{align}
\Psi(t,x,y)={}&s(t)F(x)C(y)
+a(t)\bigl(F(x)+F'(x)\bigr)D(y)
+\bigl[a(t)\bigl(F(x)+F'(x)\bigr)+s(t)F''(x)\bigr]H(y),
\label{eq:newMMS_stream}
\end{align}
we can prescribe the following exact fields
\begin{equation}
\left\{
\begin{aligned}
 \bu(t,x,y)&=\bigl(\partial_y\Psi,-\partial_x\Psi\bigr)^T,\\
 p_f(t,x,y)&=b(t)F(x)-2a(t)\bigl(F'(x)+F''(x)\bigr),\\
 \be(t,x,y)&=a(t)\bigl(F(x)A(y),F(x)B(y)\bigr)^T,\\
 \bxi(t,x,y)&=a'(t)\bigl(F(x)A(y),F(x)B(y)\bigr)^T,\\
 \phi(t,x,y)&=b(t)F(x)B(y)
 +\bigl[a'(t)F(x)+s(t)F'(x)\bigr]A(y).
\end{aligned}
\right.
\label{eq:newMMS_fields}
\end{equation}
In particularly, the stream-function representation gives fluid velocity divergence free $\nabla\cdot\bu=0$, and the
kinematic relation $\partial_t\be=\bxi$ holds identically. The factors
in~\eqref{eq:newMMS_definitions} also yield the following homogeneous exterior
conditions:
\[
\bu=0\quad\text{on }\Sigma_f,\qquad
\be=\bxi=0,\quad\phi=0\quad\text{on }\Sigma_p.
\]
At the interface $\Gamma$, these fields are nonzero and satisfy
\begin{align}
&\bu|_\Gamma
 =\bigl(a(t)(F+F'),-s(t)F'\bigr)^T,\qquad
\bxi|_\Gamma=\bigl(0,a'(t)F\bigr)^T, \qquad
\be|_\Gamma=\bigl(0,a(t)F\bigr)^T,\notag\\
&\phi|_\Gamma=b(t)F,\qquad
\partial_y\phi|_\Gamma=a'(t)F+s(t)F',\qquad
\bu_p|_\Gamma\cdot\bn_p=-a'(t)F-s(t)F',
\label{eq:newMMS_traces}
\end{align}
where $F$, $F'$, and $F''$ are evaluated at $x$. Therefore, they satisfy the corresponding coupling conditions as listed in \eqref{eq:flux}-\eqref{eq:balance}:
\[
 \bu\cdot\bn_f=(\bxi+\bu_p)\cdot\bn_f,
 \qquad \bu_p=-K\nabla\phi.
\]
\begin{align}
 \bigl(\bs_f(\bu,p_f)\bigr)_{xy}\big|_\Gamma
 &=\bigl(\bs_p(\be,\phi)\bigr)_{xy}\big|_\Gamma
 =a(t)(F+F'),\notag\\
 \bigl(\bs_f(\bu,p_f)\bigr)_{yy}\big|_\Gamma
 &=\bigl(\bs_p(\be,\phi)\bigr)_{yy}\big|_\Gamma
 =-b(t)F=-\phi|_\Gamma.
 \label{eq:newMMS_stresses}
\end{align}
In particular,
$\bxi\cdot\bn_p=a'(t)F(x)\not\equiv0$, and at the exact time levels
\begin{equation}
 \bigl(\bxi^{n+1}-2\bxi^n+\bxi^{n-1}\bigr)\cdot\bn_p
 =24t_n(\Delta t)^2F(x).
 \label{eq:newMMS_normal_second_difference}
\end{equation}
Although the exact fields satisfy the continuous physical interface conditions, their substitution into the time-discrete Robin–Robin mechanical condition produces the nonzero second-order extrapolation residual in \eqref{eq:newMMS_normal_second_difference}. 
Therefore the normal structure-velocity extrapolation term retained by the
Robin--Robin algorithm in \cite{WangLeeCanic2026SecondOrder} is nonzero for this benchmark. Moreover, the
second-order backward-difference approximation is not exact on these
higher-degree temporal polynomials. For example, when applied to the
exact structure velocity, the backward differentiation defect will be
\[
 D_t\bxi(t_{n+1})-\partial_t\bxi(t_{n+1})
 =-8(\Delta t)^2\bigl(F(x)A(y),F(x)B(y)\bigr)^T.
\]
Thus this example tests both the interface extrapolation and backward
differentiation errors, comparing to \ref{test1}.

The volume sources are obtained by substituting the exact fields into the
continuous equations:
\[
\boldsymbol f_f=\partial_t\bu-\Delta\bu+\nabla p_f,\qquad
\boldsymbol f_p=\partial_t\bxi-\Delta\be
 -2\nabla(\nabla\cdot\be)+\nabla\phi,\qquad
g_p=\partial_t\phi+\nabla\cdot\bxi-\Delta\phi.
\]
The initial data are
\begin{align*}
\bu(0)&=0,\quad p_f(0,x,y)=F(x),\quad \be(0)=0,\\
\bxi(0,x,y)&=\bigl(F(x)A(y),F(x)B(y)\bigr)^T,\qquad
\phi(0,x,y)=F(x)\bigl(A(y)+B(y)\bigr).
\end{align*}
The two initialization levels, namely $t_0=0$ and $t_1=\Delta t$,
are prescribed by projecting the exact fields. 
Same as in \ref{test1}, in both refinement studies, we use continuous piecewise quadratic ($P_2$) elements for fluid velocity, poroelastic displacement, structure velocity, and pore pressure, and continuous piecewise linear ($P_1$) elements for fluid pressure. Both subdomains are discretized using uniform triangular meshes with $n$ intervals in each coordinate direction.

\begin{table}[htbp]
\centering
\caption{Spatial convergence of the Neumann--Robin (NR, $\beta=0$) and Robin--Robin (RR, $\beta=1$) schemes for $T=1$ and fixed $\Delta t=10^{-3}$. The displacement error $E_\eta$ is measured in the $H^1$-equivalent elastic energy norm, and all other errors are measured in $L^2$. All errors are maxima over the computed BDF2 time levels.}
\label{tab:spatial-nr-rr-comparison}

\medskip
\textbf{(a) Fluid variables}\par\smallskip

\begingroup
\scriptsize
\renewcommand{\arraystretch}{1.1}
\setlength{\tabcolsep}{2pt}
\begin{tabular*}{\linewidth}{@{\extracolsep{\fill}}ccrrrrrrrr@{}}
\toprule
$n$ & $\Delta x$ & \multicolumn{4}{c}{$E_u$}
& \multicolumn{4}{c}{$E_{p_f}$}\\
\cmidrule(lr){3-6}
\cmidrule(lr){7-10}
 & & \multicolumn{2}{c}{NR} & \multicolumn{2}{c}{RR}
 & \multicolumn{2}{c}{NR} & \multicolumn{2}{c}{RR}\\
 & & Error & Rate & Error & Rate
 & Error & Rate & Error & Rate\\
\midrule
4 & $1/4$
& $2.25\times10^{0}$ & --
& $2.25\times10^{0}$ & --
& $2.33\times10^{1}$ & --
& $2.33\times10^{1}$ & --\\

8 & $1/8$
& $3.43\times10^{-1}$ & 2.71
& $3.43\times10^{-1}$ & 2.71
& $5.70\times10^{0}$ & 2.03
& $5.70\times10^{0}$ & 2.03\\

16 & $1/16$
& $4.20\times10^{-2}$ & 3.03
& $4.20\times10^{-2}$ & 3.03
& $1.31\times10^{0}$ & 2.12
& $1.31\times10^{0}$ & 2.12\\

32 & $1/32$
& $5.13\times10^{-3}$ & 3.03
& $5.13\times10^{-3}$ & 3.03
& $3.09\times10^{-1}$ & 2.08
& $3.09\times10^{-1}$ & 2.08\\

64 & $1/64$
& $6.36\times10^{-4}$ & 3.01
& $6.36\times10^{-4}$ & 3.01
& $7.56\times10^{-2}$ & 2.03
& $7.56\times10^{-2}$ & 2.03\\

128 & $1/128$
& $7.93\times10^{-5}$ & 3.00
& $7.93\times10^{-5}$ & 3.00
& $1.88\times10^{-2}$ & 2.01
& $1.88\times10^{-2}$ & 2.01\\
\bottomrule
\end{tabular*}
\endgroup

\medskip
\textbf{(b) Poroelastic variables}\par\smallskip

\begingroup
\scriptsize
\renewcommand{\arraystretch}{1.1}
\setlength{\tabcolsep}{2pt}
\begin{tabular*}{\linewidth}{@{\extracolsep{\fill}}ccrrrrrrrrrrrr@{}}
\toprule
$n$ & $\Delta x$
& \multicolumn{4}{c}{$E_\eta$}
& \multicolumn{4}{c}{$E_\xi$}
& \multicolumn{4}{c}{$E_\phi$}\\
\cmidrule(lr){3-6}
\cmidrule(lr){7-10}
\cmidrule(lr){11-14}
 & & \multicolumn{2}{c}{NR} & \multicolumn{2}{c}{RR}
 & \multicolumn{2}{c}{NR} & \multicolumn{2}{c}{RR}
 & \multicolumn{2}{c}{NR} & \multicolumn{2}{c}{RR}\\
 & & Error & Rate & Error & Rate
 & Error & Rate & Error & Rate
 & Error & Rate & Error & Rate\\
\midrule
4 & $1/4$
& $6.29\times10^{-1}$ & --
& $6.29\times10^{-1}$ & --
& $7.13\times10^{-2}$ & --
& $7.13\times10^{-2}$ & --
& $1.16\times10^{-1}$ & --
& $1.16\times10^{-1}$ & --\\

8 & $1/8$
& $1.63\times10^{-1}$ & 1.95
& $1.63\times10^{-1}$ & 1.95
& $1.25\times10^{-2}$ & 2.51
& $1.25\times10^{-2}$ & 2.51
& $1.17\times10^{-2}$ & 3.31
& $1.17\times10^{-2}$ & 3.31\\

16 & $1/16$
& $4.00\times10^{-2}$ & 2.03
& $4.00\times10^{-2}$ & 2.03
& $1.30\times10^{-3}$ & 3.27
& $1.30\times10^{-3}$ & 3.27
& $1.34\times10^{-3}$ & 3.12
& $1.34\times10^{-3}$ & 3.12\\

32 & $1/32$
& $1.00\times10^{-2}$ & 2.00
& $1.00\times10^{-2}$ & 2.00
& $1.24\times10^{-4}$ & 3.39
& $1.24\times10^{-4}$ & 3.39
& $1.66\times10^{-4}$ & 3.01
& $1.66\times10^{-4}$ & 3.01\\

64 & $1/64$
& $2.51\times10^{-3}$ & 2.00
& $2.51\times10^{-3}$ & 2.00
& $1.38\times10^{-5}$ & 3.16
& $1.42\times10^{-5}$ & 3.13
& $2.08\times10^{-5}$ & 3.00
& $2.08\times10^{-5}$ & 3.00\\

128 & $1/128$
& $6.29\times10^{-4}$ & 2.00
& $6.29\times10^{-4}$ & 2.00
& $2.86\times10^{-6}$ & 2.28
& $4.25\times10^{-6}$ & 1.74
& $2.77\times10^{-6}$ & 2.91
& $2.83\times10^{-6}$ & 2.88\\
\bottomrule
\end{tabular*}
\endgroup
\end{table}
Table~\ref{tab:spatial-nr-rr-comparison} shows approximately
second-order spatial convergence of displacement in the elastic energy norm and fluid pressure in $L^2$ for both schemes, and approximately third-order $L^2$ convergence of fluid velocity. Structure velocity and pore pressure
exhibit rates near three over much of the refinement range. On the finest refinement, the structure-velocity rates decrease to $2.28$ for NR and $1.74$ for RR, while the pore-pressure rates are $2.91$ and $2.88$, respectively. As in the preceding manufactured-solution test, the approximately third-order spatial convergence observed for several variables in the $L^2$ norm exceeds the order guaranteed by the bulk energy-norm error theorem. These numerical observations are consistent with the higher approximation order of quadratic finite elements, but they do not constitute a proof of an optimal $L^2$ error bound for either splitting scheme.

\begin{table}[htbp]
\centering
\caption{Temporal convergence of the Neumann--Robin (NR, $\beta=0$)
and Robin--Robin (RR, $\beta=1$) schemes for $T=1$ on the fixed
mesh $n=250$, $\Delta x=1/250$. The displacement error $E_\eta$
is measured in the $H^1$-equivalent elastic energy norm, and
all other errors are measured in $L^2$. All errors are maxima over the computed BDF2 time levels.}
\label{tab:temporal-nr-rr-comparison}

\medskip
\textbf{(a) Fluid variables}\par\smallskip

\begingroup
\scriptsize
\renewcommand{\arraystretch}{1.1}
\setlength{\tabcolsep}{2pt}
\begin{tabular*}{\linewidth}{@{\extracolsep{\fill}}crrrrrrrr@{}}
\toprule
$\Delta t$ & \multicolumn{4}{c}{$E_u$}
& \multicolumn{4}{c}{$E_{p_f}$}\\
\cmidrule(lr){2-5} \cmidrule(lr){6-9}
 & \multicolumn{2}{c}{NR} & \multicolumn{2}{c}{RR}
 & \multicolumn{2}{c}{NR} & \multicolumn{2}{c}{RR}\\
 & Error & Rate & Error & Rate
 & Error & Rate & Error & Rate\\
\midrule
$1/4$  & $2.30\times10^{-2}$ & --
       & $2.31\times10^{-2}$ & --
       & $2.02\times10^{-1}$ & --
       & $2.10\times10^{-1}$ & --\\
$1/8$  & $6.61\times10^{-3}$ & 1.80
       & $6.65\times10^{-3}$ & 1.80
       & $6.26\times10^{-2}$ & 1.69
       & $6.72\times10^{-2}$ & 1.64\\
$1/16$ & $1.77\times10^{-3}$ & 1.90
       & $1.77\times10^{-3}$ & 1.91
       & $1.76\times10^{-2}$ & 1.83
       & $1.87\times10^{-2}$ & 1.85\\
$1/32$ & $4.56\times10^{-4}$ & 1.95
       & $4.58\times10^{-4}$ & 1.95
       & $6.60\times10^{-3}$ & 1.42
       & $6.79\times10^{-3}$ & 1.46\\
$1/64$ & $1.16\times10^{-4}$ & 1.97
       & $1.17\times10^{-4}$ & 1.97
       & $5.04\times10^{-3}$ & 0.39
       & $5.06\times10^{-3}$ & 0.42\\
\bottomrule
\end{tabular*}
\endgroup

\medskip
\textbf{(b) Poroelastic variables}\par\smallskip

\begingroup
\scriptsize
\renewcommand{\arraystretch}{1.1}
\setlength{\tabcolsep}{2pt}
\begin{tabular*}{\linewidth}{@{\extracolsep{\fill}}crrrrrrrrrrrr@{}}
\toprule
$\Delta t$ & \multicolumn{4}{c}{$E_\eta$}
& \multicolumn{4}{c}{$E_\xi$}
& \multicolumn{4}{c}{$E_\phi$}\\
\cmidrule(lr){2-5}
\cmidrule(lr){6-9}
\cmidrule(lr){10-13}
 & \multicolumn{2}{c}{NR} & \multicolumn{2}{c}{RR}
 & \multicolumn{2}{c}{NR} & \multicolumn{2}{c}{RR}
 & \multicolumn{2}{c}{NR} & \multicolumn{2}{c}{RR}\\
 & Error & Rate & Error & Rate
 & Error & Rate & Error & Rate
 & Error & Rate & Error & Rate\\
\midrule
$1/4$  & $3.06\times10^{-1}$ & --
       & $2.85\times10^{-1}$ & --
       & $8.63\times10^{-2}$ & --
       & $1.73\times10^{-1}$ & --
       & $3.82\times10^{-2}$ & --
       & $4.57\times10^{-2}$ & --\\
$1/8$  & $1.04\times10^{-1}$ & 1.56
       & $9.66\times10^{-2}$ & 1.56
       & $2.82\times10^{-2}$ & 1.61
       & $5.53\times10^{-2}$ & 1.64
       & $1.20\times10^{-2}$ & 1.67
       & $1.48\times10^{-2}$ & 1.63\\
$1/16$ & $2.97\times10^{-2}$ & 1.81
       & $2.76\times10^{-2}$ & 1.81
       & $8.13\times10^{-3}$ & 1.80
       & $1.48\times10^{-2}$ & 1.90
       & $3.31\times10^{-3}$ & 1.86
       & $4.01\times10^{-3}$ & 1.88\\
$1/32$ & $7.90\times10^{-3}$ & 1.91
       & $7.39\times10^{-3}$ & 1.90
       & $2.16\times10^{-3}$ & 1.91
       & $3.77\times10^{-3}$ & 1.97
       & $8.64\times10^{-4}$ & 1.94
       & $1.04\times10^{-3}$ & 1.95\\
$1/64$ & $2.04\times10^{-3}$ & 1.95
       & $1.92\times10^{-3}$ & 1.95
       & $5.54\times10^{-4}$ & 1.96
       & $9.49\times10^{-4}$ & 1.99
       & $2.20\times10^{-4}$ & 1.97
       & $2.63\times10^{-4}$ & 1.98\\
\bottomrule
\end{tabular*}
\endgroup
\end{table}
Table~\ref{tab:temporal-nr-rr-comparison} reports the temporal convergence results on the fixed mesh $n=250$, with $T=1$ and $\Delta t=1/4,1/8,1/16,1/32,1/64$. For both schemes, the convergence rates for displacement, structure velocity, and pore pressure approach two as the time step decreases, supporting the second-order temporal accuracy of the BDF2--EX2 discretization for these variables. On the finest refinement, the observed rates for $(E_\eta,E_\xi,E_\phi)$ are $(1.95,1.96,1.97)$ for the NR scheme and $(1.95,1.99,1.98)$ for the RR scheme, respectively. The fluid-velocity errors also exhibit nearly second-order convergence, with final observed rates of $1.97$ for both schemes. In contrast, the fluid-pressure errors decrease initially but gradually level off near $5.04\times10^{-3}$ for the NR scheme and $5.06\times10^{-3}$ for the RR scheme. The final observed convergence rates decrease to $0.39$ and $0.42$, respectively. This reduction in the observed temporal convergence rate of fluid pressure at smaller time steps is consistent with the influence of the fixed spatial discretization error, which becomes dominant as the temporal discretization error decreases.

In summary, for identical spatial meshes and time steps, the NR scheme consistently yields smaller errors in the structure velocity and pore pressure across all temporal refinement cases, while the RR scheme achieves smaller displacement errors. The fluid-velocity errors are comparable between the two methods. Overall, neither scheme demonstrates a uniform accuracy advantage across all solution variables.

\subsection{Numerical stability study based on 3D cylindrical pressure-pulse benchmark}
\label{sec:three_dimensional_pressure_pulse}
We want to investigate the stability of the NR and
RR schemes by fixing the mesh and increasing the time step. The theoretical sufficient restrictions
have the forms as follows:
\[
\Delta t \leq C_{\mathrm{NR}}h,
\qquad
\Delta t \leq C_{\mathrm{RR}}h^2.
\]

We compare the monolithic, second-order NR, and
RR methods in a 3D setting by considering a
cylindrical fluid lumen surrounded by a poroelastic wall and subjected to a inlet-pressure pulse. All three methods use the same reference mesh, physical parameters, loading, and exterior boundary conditions. The wall is modeled using small-deformation poroelasticity on a fixed reference geometry: its displacement is computed, but the fluid mesh is not moved.
\begin{figure}[htbp]
    \centering
    \includegraphics[width=0.95\textwidth]{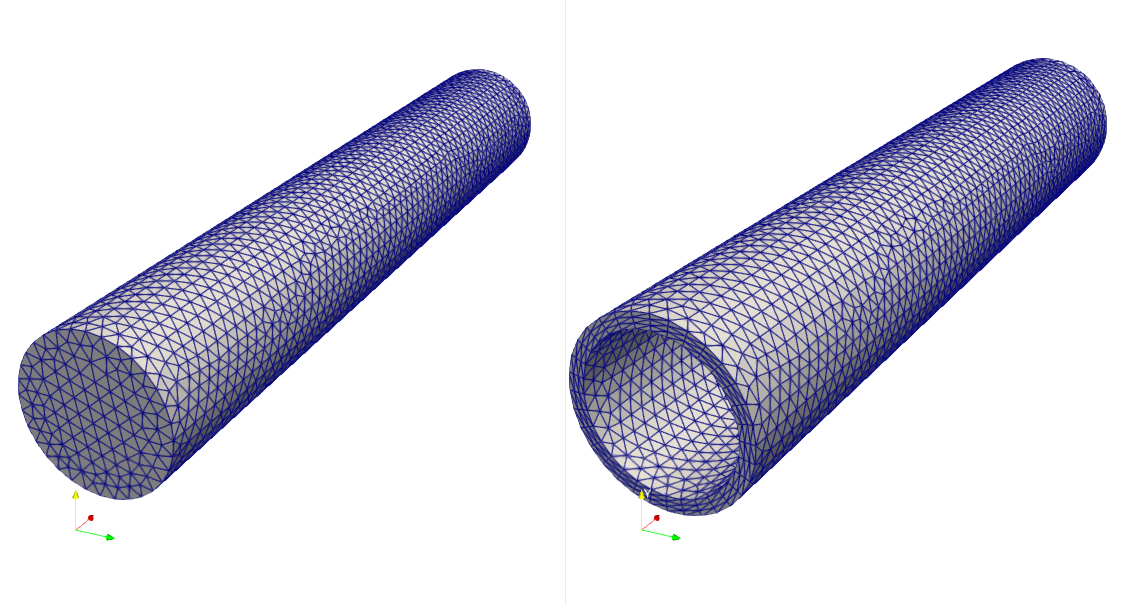}
    \caption{3D computational meshes for the
    cylindrical Stokes--Biot problem.
    Left: tetrahedral mesh of the fluid lumen
    $\Omega_f$.
    Right: tetrahedral mesh of the surrounding
    poroelastic wall $\Omega_p$, consisting of
    three radial layers. 
    The fluid and structure meshes share a
    conforming interface at $r=0.5\,\mathrm{cm}$. The fluid mesh contains $20{,}892$ tetrahedra and the wall mesh contains
$40{,}518$ tetrahedra, giving $61{,}410$ tetrahedra and $11{,}668$
unique mesh vertices in total. The interface has $4{,}502$ matching
triangular facets.}
    \label{fig:pulse_mesh}
\end{figure}

The cylinder axis is aligned with the $x$-direction and
$r=\sqrt{y^2+z^2}$. The fluid and poroelastic domains are
\begin{equation}
\begin{aligned}
\Omega_f=\bigl\{(x,y,z):0<x<6,\;r<0.5\bigr\},\qquad
\Omega_p=\bigl\{(x,y,z):0<x<6,\;0.5<r<0.6\bigr\},
\end{aligned}
\label{eq:pulse_domains}
\end{equation}
where all lengths are in centimeters. Their common interface is
$\Gamma=\{(x,y,z):0<x<6,\;r=0.5\}$; the inlet and outlet are the
cross-sections at $x=0$ and $x=6$, respectively. 
% \subsection{Governing equations and material parameters}
% \label{subsec:pulse_equations}
% The fluid satisfies the time-dependent incompressible Stokes system
% \begin{equation}
% \rho_f\partial_t\boldsymbol u-\nabla\cdot\boldsymbol\sigma_f=\boldsymbol0,
% \qquad \nabla\cdot\boldsymbol u=0\quad\text{in }\Omega_f,
% \label{eq:pulse_stokes}
% \end{equation}
% with
% \[
% \boldsymbol\sigma_f=2\mu_f\boldsymbol D(\boldsymbol u)-p_f\boldsymbol I,
% \qquad
% \boldsymbol D(\boldsymbol v)=\tfrac12\bigl(\nabla\boldsymbol v
% +\nabla\boldsymbol v^{T}\bigr).
% \]
% The poroelastic wall satisfies the dynamic Biot system
% \begin{equation}
% \begin{aligned}
% \partial_t\boldsymbol\eta&=\boldsymbol\xi,\\
% \rho_p\partial_t\boldsymbol\xi-\nabla\cdot\boldsymbol\sigma_p
% &=\boldsymbol0,\\
% C_0\partial_t\phi+\alpha\nabla\cdot\boldsymbol\xi
% -\nabla\cdot(K\nabla\phi)&=0
% \end{aligned}
% \qquad\text{in }\Omega_p,
% \label{eq:pulse_biot}
% \end{equation}
% where
% \[
% \boldsymbol\sigma_p
% =2\mu_p\boldsymbol D(\boldsymbol\eta)
% +\lambda_p(\nabla\cdot\boldsymbol\eta)\boldsymbol I
% -\alpha\phi\boldsymbol I.
% \]
% Here $\boldsymbol\eta$, $\boldsymbol\xi$, and $\phi$ denote wall
% displacement, wall velocity, and pore pressure, respectively. Darcy flux is
% $\boldsymbol q=-K\nabla\phi$; the hydraulic mobility $K$ already includes
% the appropriate viscosity factor and is not divided by $\mu_f$ again.
All quantities are expressed in cgs units. The
parameters used in the simulations are given in Table~\ref{tab:pulse_parameters}.
These values are taken from the numerical model and are not presented as a
calibration to experimental vessel data.
\begin{table}[htbp]
\centering
\caption{Physical and numerical parameters of the cylindrical
pressure-pulse benchmark.}
\label{tab:pulse_parameters}
\small
\begin{tabular}{lll}
\toprule
Quantity & Value & Units \\
\midrule
Fluid density $\rho_f$ & $1.0$ & $\mathrm{g}/\mathrm{cm}^3$ \\
Wall density $\rho_p$ & $1.1$ & $\mathrm{g}/\mathrm{cm}^3$ \\
Fluid viscosity $\mu_f$ & $0.035$ & $\mathrm{g}/(\mathrm{cm}\,\mathrm{s})$ \\
Wall shear modulus $\mu_p$ & $1.6779\times10^6$ & $\mathrm{dyn}/\mathrm{cm}^2$ \\
Wall Lam\'e coefficient $\lambda_p$ & $8.2215\times10^7$ & $\mathrm{dyn}/\mathrm{cm}^2$ \\
Biot coefficient $\alpha$ & $1$ & dimensionless \\
Storage coefficient $C_0$ & $10^{-3}$ & $\mathrm{cm}^2/\mathrm{dyn}$ \\
Isotropic Darcy mobility $K$ & $10^{-8}\boldsymbol I$ & $\mathrm{cm}^3\,\mathrm{s}/\mathrm{g}$ \\
Tangential coefficient $\gamma$ & $1000$ & $\mathrm{dyn}\,\mathrm{s}/\mathrm{cm}^3$ \\
Numerical Robin parameter $L$ & $1500$ & $\mathrm{dyn}\,\mathrm{s}/\mathrm{cm}^3$ \\
\bottomrule
\end{tabular}
\end{table}

A spatially uniform inlet pressure is applied at the inlet
$\Gamma_{\mathrm{in}}=\{x=0,r<0.5\}$:
\begin{equation}
p_{\mathrm{in}}(t)=
\begin{cases}
\dfrac{13333}{2}\left[1-\cos\left(\dfrac{2\pi t}{0.003}\right)\right],
&0\le t\le 0.003\,\mathrm{s},\\[5pt]
0,&t>0.003\,\mathrm{s}.
\end{cases}
\label{eq:pulse_inlet_pressure}
\end{equation}
and we have the following boundary conditions:
\begin{equation}
\begin{aligned}
\boldsymbol\sigma_f\boldsymbol n_f
  &=-p_{\mathrm{in}}(t)\boldsymbol n_f
   &&\text{on the fluid inlet},\\
\boldsymbol\sigma_f\boldsymbol n_f
  &=\boldsymbol0 &&\text{on the fluid outlet},\\
\boldsymbol\eta=\boldsymbol\xi
  &=\boldsymbol0 &&\text{on the two wall-end annuli},\\
\boldsymbol\sigma_p\boldsymbol n_p
  &=\boldsymbol0 &&\text{on the outer wall }r=0.6,\\
\boldsymbol q\cdot\boldsymbol n_p
  &=0 &&\text{on all exterior porous boundaries}.
\end{aligned}
\label{eq:pulse_exterior_conditions}
\end{equation}

\begin{figure}[htbp]
    \centering

    \begin{subfigure}[t]{0.42\textwidth}
        \centering
        \includegraphics[width=\linewidth]{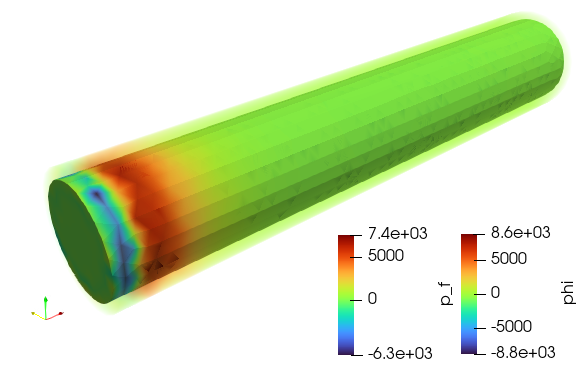}
        \caption{$t=6~\mathrm{ms}$}
        \label{fig:pressure_6ms}
    \end{subfigure}
 \quad
    \begin{subfigure}[t]{0.42\textwidth}
        \centering
        \includegraphics[width=\linewidth]{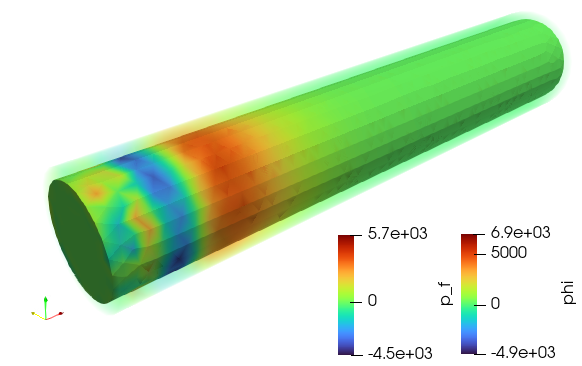}
        \caption{$t=12~\mathrm{ms}$}
        \label{fig:pressure_12ms}
    \end{subfigure}

    \vspace{0.2em}

    \begin{subfigure}[t]{0.42\textwidth}
        \centering
        \includegraphics[width=\linewidth]{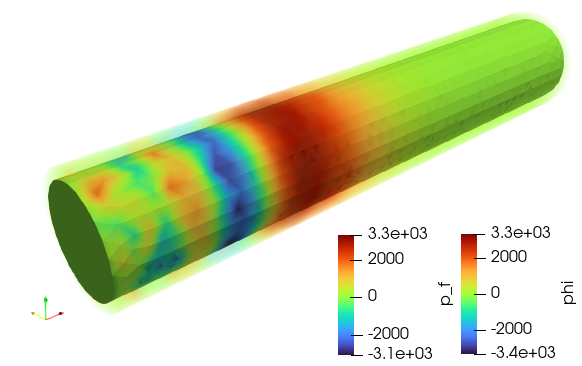}
        \caption{$t=18~\mathrm{ms}$}
        \label{fig:pressure_18ms}
    \end{subfigure}
\quad
    \begin{subfigure}[t]{0.42\textwidth}
        \centering
        \includegraphics[width=\linewidth]{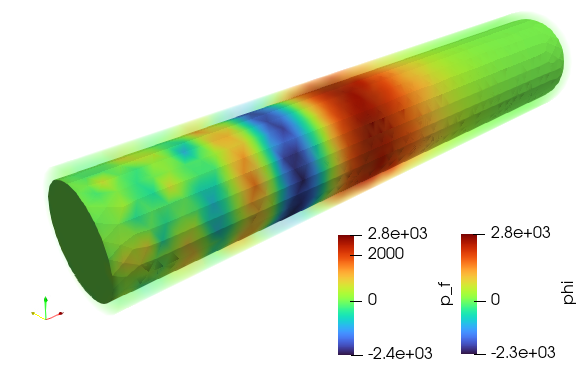}
        \caption{$t=24~\mathrm{ms}$}
        \label{fig:pressure_24ms}
    \end{subfigure}

    \caption{
    Snapshots of the fluid pressure $p_f$
    and pore pressure $\phi$ at
    (a) $t=6~\mathrm{ms}$,
    (b) $t=12~\mathrm{ms}$,
    (c) $t=18~\mathrm{ms}$, and
    (d) $t=24~\mathrm{ms}$, obtained from the NR scheme.
    The pressure fields are overlaid to illustrate
    the coupled response of the fluid lumen and
    surrounding poroelastic wall.
    Color scales vary between snapshots.
    }
    \label{fig:pressure_snapshots}
\end{figure}

Figure~\ref{fig:pressure_snapshots} shows the overlaid
pressure fields at $t=6$, $12$, $18$, and $24~\mathrm{ms}$, obtained from the NR scheme. 
The fluid pressure is displayed in the cylindrical lumen, while the
pore pressure is shown in the surrounding poroelastic
wall.

\begin{remark}[Scope of the stability analysis]
The 3D benchmark employs boundary conditions
that differ from the homogeneous exterior Dirichlet conditions
assumed in the theoretical stability analysis. In particular,
the benchmark includes nonhomogeneous fluid traction loading,
a traction-free mechanical boundary, and an impermeable
exterior porous boundary. Consequently, the stability theorem
does not directly apply to this configuration. The numerical
experiment is intended to assess the behavior of the proposed
Neumann-Robin splitting under these more general boundary
conditions, rather than to provide a direct verification of
the theoretical time-step restriction.
\end{remark}

\paragraph{Numerical stability under time-step refinement.}

We investigate the stability of the NR
and RR splitting schemes by comparing
their physical energy evolution for different time-step
sizes on the same 3D cylindrical mesh.
All simulations use continuous $P_2$ finite elements
for the fluid velocity $\boldsymbol u$, displacement
$\boldsymbol\eta$, structure velocity $\boldsymbol\xi$,
and pore pressure $\phi$, and continuous $P_1$
elements for the fluid pressure $p_f$.
The time steps are
\begin{equation}
\Delta t\in
\left\{
1\times10^{-5},\,
2\times10^{-5},\,
2.5\times10^{-5},\,
5\times10^{-5}
\right\}\,\mathrm{s},
\qquad
T=0.024\,\mathrm{s}.
\label{eq:pulse_time_settings}
\end{equation}
Both splitting schemes use identical spatial meshes,
physical parameters, loading, and boundary conditions.
All initial fields and volume sources are zero.
The monolithic solution computed with
$\Delta t=10^{-5}\,\mathrm{s}$ is included as a
reference.

Since the inlet-pressure pulse supplies mechanical
work during $0<t<0.003\,\mathrm{s}$, an increase in
energy within this interval does not necessarily
indicate numerical instability. The instantaneous
inlet power is given by
\begin{equation}
P_{\mathrm{in}}(t)
=-\int_{\Gamma_{\mathrm{in}}}
p_{\mathrm{in}}(t)\,
\boldsymbol u\cdot\boldsymbol n_f\,\mathrm dS.
\label{eq:pulse_inlet_power}
\end{equation}
We therefore focus on the energy behavior after
the pulse and compare the splitting solutions
with the monolithic reference.

\begin{figure}[htbp]
    \centering
    \includegraphics[width=0.9\textwidth]
    {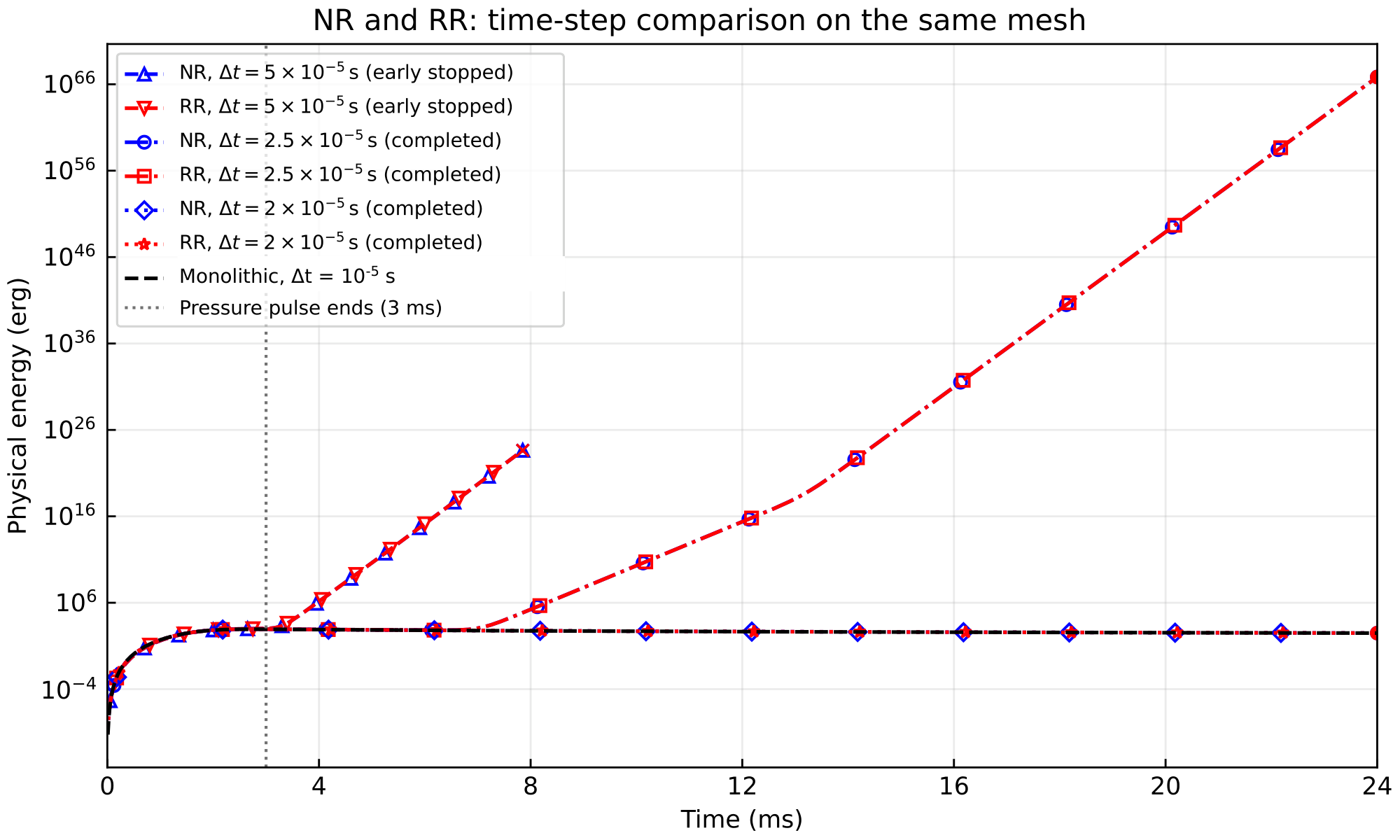}
    \caption{
    Physical energy evolution of the NR and RR schemes for different
    time-step sizes on the same 3D
    cylindrical mesh. The dashed black curve
    represents the small-step monolithic reference
    solution, and the vertical dotted line marks
    the end of the inlet-pressure pulse at
    $t=3\,\mathrm{ms}$.
    The two splitting schemes exhibit nearly
    identical energy behavior for the tested
    time steps.
    }
    \label{fig:NR_RR_stability_comparison}
\end{figure}

As shown in Figure~\ref{fig:NR_RR_stability_comparison},
both NR and RR exhibit rapid energy growth for
$\Delta t=5\times10^{-5}\,\mathrm{s}$, leading to
early termination. For
$\Delta t=2.5\times10^{-5}\,\mathrm{s}$, both schemes
complete the simulations but develop substantial
nonphysical energy amplification after the
inlet-pressure pulse. In contrast, the solution obtained from a smaller step of $\Delta t=2\times10^{-5}\,\mathrm{s}$ 
remains bounded throughout
the simulation. The NR and RR energy curves
nearly coincide at each tested time step,
indicating no discernible difference in their
numerical stability behavior on the present mesh.

These findings highlight the distinction between
theoretical sufficient stability conditions and
actual numerical stability thresholds. The NR
analysis establishes stability under a linear
restriction $\Delta t=O(h)$, improving upon the
parabolic restriction $\Delta t=O(h^2)$ obtained
in the previous RR analysis. However, neither
condition is known to be necessary or sharp.
The nearly identical numerical behavior suggests
that the more restrictive theoretical condition
for RR may be conservative and does not necessarily
reflect its practical stability limit. Thus, although NR admits a less restrictive rigorous stability guarantee, the present
experiments do not demonstrate a corresponding
numerical stability advantage over RR.

% For each scheme, we monitor the normalized bulk BDF2 energy
% and its maximum and final amplification:
% \[
% R^j=\frac{\mathcal E^j}{\mathcal E^1},\qquad
% G_{\max}=\max_{1\leq j\leq N}R^j,\qquad
% G_{\mathrm{final}}=R^N.
% \]
% The four energy contributions are recorded separately.
% Sustained amplification, nonfinite values, and excessive
% growth are examined; an isolated moderate energy increase
% is not sufficient to declare instability.

\section{Conclusion}
\label{sec:conclusion}
We have formulated and analyzed a Neumann--Robin modification of the second-order Robin-based splitting in \cite{WangLeeCanic2026SecondOrder}. The fluid normal Robin condition, Darcy Robin condition, and tangential Robin/BJS treatment are retained, while the Biot mechanical normal condition is imposed with the current pore pressure and without the auxiliary normal structure-velocity Robin term. The fluid and Biot subproblems remain independent at each time step and can be solved in parallel using extrapolated interface data. 

The removal of the extrapolated normal structure-velocity residual permits a fully discrete energy stability estimate under the linear time-step restriction $\Delta t\le C_*h$ for fixed coefficients satisfying the stated assumptions. The stability result covers both homogeneous and nonhomogeneous source terms. Under appropriate regularity, finite element approximation, and initialization assumptions, we also establish a bulk energy-norm error estimate of order $h^k+\Delta t^2$ for the fluid velocity, structure velocity, pore pressure, and displacement. The error analysis employs a common elasticity projection and a BDF2 summation-by-parts argument to control the pore-pressure projection traction term without requiring normal structure-velocity dissipation.

The manufactured-solution experiments exhibit the expected spatial convergence behavior and temporal convergence trends consistent with second-order accuracy. Reduced convergence rates on the finest temporal refinements are partly attributable to spatial-error saturation. In the 3D cylindrical benchmark, the Neumann-Robin and Robin-Robin schemes exhibit nearly identical numerical stability behavior for the tested time steps. The improvement from a parabolic to a linear mesh-dependent sufficient stability condition is therefore established analytically, while the numerical experiments do not demonstrate a corresponding stability advantage or determine the sharp stability threshold of either scheme.

\section*{Acknowledgments}
Wang's research has been supported in part by the National Science Foundation under grant DMS-2247001 and by a Simons Foundation Travel Award.

\bibliographystyle{unsrt}
\bibliography{references}
\end{document}